\documentclass{gtpart}   
\usepackage{geometry}                
\usepackage{graphicx}
\usepackage{amssymb}
\usepackage{epstopdf}
\usepackage{wrapfig }
\usepackage{cutwin}

\usepackage{caption}
\usepackage{amsmath}
\usepackage{amsthm}
\usepackage{delarray}

\usepackage[all]{xy}
\usepackage{tikz}

\usepackage{pinlabel}
\usepackage[all]{xy}

\title{Exotic symplectic manifolds from Lefschetz fibrations.}

\author{Maksim Maydanskiy} 
\givenname{Maksim}
\surname{Maydanskiy}
\address{Department of Mathematics\\Stanford University\\\newline
         Stanford, CA 94305\\USA}
\email{maksimm@math.stanford.edu}
\urladdr{http://math.stanford.edu/~maksimm/}

\keyword{example}
\keyword{sample layout}
\subject{primary}{msc2010}{53D40}
\subject{secondary}{msc2010}{32Q28}
\subject{secondary}{msc2010}{57R17}

\arxivreference{}  
\arxivpassword{}   

\newtheorem{thm}{Theorem}[section]   
\newtheorem{lem}[thm]{Lemma}         
 
\theoremstyle{definition}
   
\newtheorem{rem}[thm]{Remark}

\makeop{Homo}
\numberwithin{equation}{section}

\newtheorem*{tm}{Theorem}

\newtheorem{prop}[thm]{Proposition}

\theoremstyle{remark}
\newtheorem{lemma}[thm]{Lemma}
\newtheorem{obs}[thm]{Observation}

\theoremstyle{definition}
\newtheorem{dfn}[thm]{Definition}
\newtheorem{ex}[thm]{Example}

\newcommand{\CC}{{\mathbb C}}

\usepackage{titlesec}

\titleformat{name=\section}[hang]{\center\sc}{\thesection. }{0pt}{}[]
\titleformat{name=\subsection}[hang]{\raggedright \bfseries}{\thesubsection.}{0pt}{ }[]

\begin{document}

\begin{abstract}
In this paper we construct, in all odd complex dimensions, pairs of Liouville domains
$W_0$ and $W_1$ which are diffeomorphic to the cotangent bundle of the sphere
with one extra subcritical handle, but are not symplectomorphic. While $W_0$ is symplectically very similar to the cotangent
bundle itself, $W_1$ is more unusual. We use Seidel's exact triangles for Floer cohomology to show that the wrapped Fukaya category of $W_1$ is trivial. As a corollary we obtain that $W_1$ contains no
compact exact Lagrangian submanifolds.
\end{abstract}

\maketitle

\begin{section}{Introduction.}

This paper is concerned with the symplectic topology properties of Liouville domains, also known, from the perspective of complex geometry, as Stein domains. We use Lefschetz fibrations to construct such manifolds and wrapped Floer cohomology to distinguish their symplectomorphism types.

Lefschetz pencils were introduced by Donaldson \cite{DICM} and shortly afterwards Lefschetz theory emerged as a significant tool in symplectic topology, as manifest in the work of Auroux \cite{BRCOV}, 
Akbulut-Ozbagci \cite{AO}, Gompf \cite{Gompf}, \cite{G}, Seidel \cite{S2}, \cite{Book} and others. An important feature of Lefschetz fibrations is that they encode the topology of the total space in terms of the fiber and a collection of submanifolds in it,  the so called vanishing cycles. In what may be called an opposite direction, this allows one to construct symplectic manifolds from collections of vanishing cycles in the fiber. In particular,  modifying a given collection one gets various families of total spaces. This paper is concerned with understanding symplectic invariants of the resulting manifolds, in particular their wrapped Floer cohomology and wrapped Fukaya category. We prove the following theorem.

\begin{thm} 
\label{mainT}
There exist Liouville domains $W_0$ and $W_1$ where $W_i$ is obtained by attaching an $n$-handle to $T^{*} S^{n+1}$ ($n$ even, $n \geq 2$), with  the following properties: 
\begin{itemize}
\item $W_0$ and $W_1$ are diffeomorphic
\item $W_0$ and $W_1$ carry Lefschetz fibrations over the disc, such that the wrapped Floer cohomologies of the Lefschetz thimbles are non-zero for $W_0$ and zero for $W_1$.
\item $W_1$ does not contain any compact exact Lagrangian submanifolds (whereas $W_0$ contains a Lagrangian $S^{n+1}$, which is exact since $S^{n+1}$ is simply connected)
\end{itemize}
\end{thm}

In particular $W_0$ and $W_1$ are not exact deformation equivalent; $W_1$ is the ``exotic'' version of $W_0$. After this work has been completed, Abouzzaid and Seidel in \cite{AS2} proved that such doppelg\"angers with vanishing wrapped Fukaya categories exist for all total spaces of Lefschetz fibrations of real dimension 6 and above.

The rest of this paper is structured as follows. In section 2 we introduce Liouville domains and exact symplectic manifolds with corners and summarize their basic properties. Sections 3 and 4 are concerned with exact Lefschetz fibrations for Liouville domains. We discuss vanishing cycles, Dehn twists, thimbles, and matching cycles. In section 5 we proceed to construct our main objects of study - the Liouville domains
$W_0$ and $W_1$ - as total spaces of Lefschetz fibrations.  Section 6 reviews the basics of wrapped Floer cohomology and wrapped Fukaya categories as relevant to Lefschetz fibrations.  Section 7 recalls Seidel's exact triangles in Floer cohomology and applies them to $W_0$ and $W_1$, proving the main non-symplectomorphism result. In the final section we put this paper into a general framework of  computations of Floer-theoretic invariants of Liouville domains and discuss extensions of the present work to a more general setting.

\end{section}

\begin{section}{Liouville domains and exact symplectic manifolds with corners.}

We study the symplectic topology of Liouville domains. The introduction below closely follows \cite{Biased}. 

\begin{dfn} A \emph{Liouville domain} is a compact manifold with
boundary $M^{2n}$, together with a one-form $\theta$ which has the following two
properties. Firstly, $\omega=d \theta$ should be symplectic. Secondly, the vector
field $Z$ defined by $i_Z \omega = \theta$ should point strictly outwards along $\partial M$.
 
\end{dfn}

\begin{ex} \label{stein}

A Stein manifold $U$ with complex structure $J$ admits an exhausting  function $h: U \to \mathbb{R}$ which is strictly plurisubharmonic, meaning that $-d d_c h = - d(d h \cdot J)$ is a K\"ahler form. Then, if $C$ is a
regular value of $h$, the sublevel set $M=  h^{-1} ((-\infty; C])$ is a Liouville domain with $\theta= d_c h$ and the Liouville vector field $Z$ is the gradient of $h$ with respect to the K\"ahler metric.

\end{ex}

Note that $\alpha = \theta|_{\partial M}$ is a contact form on $\partial M$, and the negative time flow of $Z$ defines a canonical collar neighborhood  $\kappa: (- \infty, 0] \times \partial M \to M $, with $\kappa^* \theta = e^s \alpha$, $ \kappa^* Z = \partial_s$. This collar is modeled on the negative part of symplectization of $(\partial M, \alpha)$ and allows us to complete $M$ by attaching an infinite cone corresponding to the positive half:

$\hat{M}= M \cup (\partial M \times [0, \infty))$; $\hat{\theta}|([0, \infty)\times M)= e^s \alpha$; $\hat{Z}|([0, \infty)\times M)= \partial_s.$

Such a non-compact $\hat{M}$ will be called a Liouville manifold.

We will sometimes use the coordinate $r=e^s$ on the cone.

\begin{dfn}
 
A \emph{Liouville isomorphism} between domains $M_0$ and $M_1$ is a diffeomorphism
$\phi : \hat{M}_0 \to \hat{M}_1$ satisfying $\phi^*\theta_1 = \theta_0 +d g$, where $g$ is a compactly supported smooth function.  
\end{dfn}

Note that $\phi$ is defined on the completions, rather than the domains themselves. Any such $\phi$ is a symplectomorphism and is compatible with the
Liouville flow at infinity. This means that on a piece of the cone $[\rho, \infty) \times \partial M_0 \subset \hat{M}_0$ for
some $\rho > 0$, it has the form
$\phi(s, y) = (s - f(y), \psi(y))$,
where $\psi: \partial M_0 \to  \partial M_1$ is a contact isomorphism, satisfying $\psi^* \alpha_1 =
e^f \alpha_0$  for some function $f$. Thus isomorphic Liouville domains have contactomorphic boundaries. Note that the contact form on the boundary  is not preserved and in fact can be changed arbitrarily within the class defining the same contact structure.

\begin{ex}

Let $N$ be  a manifold and $T^*N$ its cotangent bundle with the standard symplectic form $\omega=d\lambda=\Sigma dp \wedge dq$. Then the vector field $Z=\Sigma p \frac{ \partial}{\partial p}$ generates the Liouville flow of ``radial rescaling''. Any choice of metric on $N$ makes the corresponding unit disc bundle into a Liouville domain. All such domains are Liouville isomorphic, with corresponding completions symplectomorphic to $T^*N$ itself.

\end{ex}

The following version of Moser's Lemma, which says that
deformation equivalence implies Liouville isomorphism, holds in this context.
\begin{lem}

 Let $(\theta_t)_{0\leq t \leq 1}$ be a smooth family of Liouville structures on M.
Then all the $(M, \theta_t)$ and $(M, \theta_t')$ are  Liouville isomorphic for any $t, t' \in [0,1]$.
\end{lem}

\begin{ex}

If in Example \ref{stein} the critical point set of $h:U \to \mathbb{R}$ is compact, then taking $C$ to
be bigger than the largest critical value, we get a Liouville domain which is independent of the particular choice of $C$ up to Liouville isomorphism.
If we assume in addition that $h$ is complete, then $(U,-d d_c h)$ itself
will be symplectically isomorphic to $\hat{M}$.
 In this context completeness of the gradient vector field can always be
achieved by a reparametrization $h \to \beta(h)$, where $\beta:\mathbb{R} \mapsto \mathbb{R}$ is a function with positive first and second derivatives, that can be constructed explicitely, see  \cite{Biran}, Lemma 3.1.

\end{ex}

We will want to do Floer theory on Liouville manifolds. Do do this one needs to control the behavior of holomorphic maps near infinity  - that is on the attached cone $\partial M \times [0, \infty)$. To achieve this one needs to impose some control on the almost-complex structure on $\hat{M}$. This motivates the following definition.

\begin{dfn} 
\label{contactJ} (\cite{Biased}, Section 3c) An almost complex structure $J$ on $\hat{M}$ is called \emph{of contact type} if $d(e^s) \cdot J=-\theta$.  \end{dfn}

The notion of Liouville domain is closely related to that of \emph{exact symplectic manifold with corners}, as defined in \cite{Book}, Section 7a.
We introduce some terminology.

A \emph{smooth manifold with corners} $M$ is covered by charts of the form $[0,\infty) \times \mathbb{R}^{n-k}$, and for a point $x\in M$ given in such a chart by coordinates $(x_1, \ldots x_n)$, one defines the \emph{depth} $d(x)$ to be the number of $(x_1, \ldots x_k)$ which are zero. This is independent of the choice of the chart. Connected components of the union of all points of depth $k$ are called the \emph{connected boundary strata} of depth $k$, and their closures in $M$ are called \emph{connected boundary faces} of depth $k$. We require that these faces have no self-intesections (thus $M$ is a manifold with faces, in the terminology of Remark 2.11 in \cite{Joyce}), so that the point of depth $k$ lies in exactly $2^k -1$ connected boundary faces (of all positive depths). Without the connectedness condition, a \emph{boundary face} is a union of pairwise disjoint connected faces of same depth; and a \emph{boundary stratum} is a union of connected boundary strata of the same depth, having disjoint closures.

\begin{dfn} \label{def-MwC} An \emph{exact symplectic manifold with corners} $(M, \omega_M, \theta_M, I_M)$ is a compact smooth
manifold with corners $M$, equipped with a symplectic form $\omega_M$, a one-form $\theta_M$
satisfying $d \theta_M = \omega_M$, and an $\omega_M$-compatible almost complex structure $I_M$. These
should satisfy two convexity conditions: the Liouville vector field must point strictly outwards along all boundary faces of $\partial_M$; and the boundary must
be weakly $I_M$-convex, which means that $I_M$-holomorphic curves cannot touch $\partial M$
unless they are completely contained in it. 
\end{dfn}

Exact symplectic manifolds are technically more convenient when working with fibrations. Note that a Liouville domain $M$ with a choice of compatible almost complex structure becomes an exact symplectic manifold (without corners) - all conditions except weak boundary convexity are automatic, and the maximum principle for holomorphic curves ensures that last condition as well. In the opposite direction, the only thing that will be important to us is that exact symplectic manifolds obtained in the course of our constructions will have at most codimension one corners and that  such corners can be smoothed to make the resulting manifolds into honest Liouville domains (this is Lemma 7.6 in \cite{Book}; a similar smoothing occurs in the process of Weinstein handle attachment, \cite{Wein}). All the invariants that we will consider will be insensitive to the details of these smoothings, so long as they happen in sufficiently small neighborhoods of the corners.

\end{section}

\begin{section}{Lefschetz fibrations.}

 Simply put, a Lefschetz fibration is a map with isolated singularities modeled on the complex singularity of the simplest type. The discussion below formalizes this description for the category of Liouville domains.

Most of the technical setup follows \cite{Book}. The sections most relevant for us are 15 and 16. We summarize what will be needed below.

We will denote by  $\mathbb{D}^2$ the standard disc in $\mathbb{R}$ of radius $\sqrt {2}$. We make it into a Liouville domain by choosing the one-form $\eta=\frac{1}{2} r_{polar}^2 d \theta$, where $(r_{polar}, \theta)$ are  the polar coordiantes, so   that $\eta_{\partial \mathbb{D}^2}> 0$ and $\omega_{st} = d \eta$ is the standard symplectic form. The completion $\hat{\mathbb{D}}^2$ is exact symplectomorphic to $\mathbb{C}$ with its standard symplectic form and primitive via the identification that extends the obvious inclusion of $\mathbb{D}^2$ to $\mathbb{C}$ and on the conical end sends $(e^s, \theta)$ to $(\frac{1}{2} r_{polar}^2,  \theta)$. Then the standard complex structure on $\mathbb{C}$ is of contact type in the sense of Definition \ref{contactJ} (this explains the choice of $\sqrt{2}$ as the radius of $\mathbb{D}^2$).

\begin{dfn}\label{def-LF} Let $E$ An (exact) \emph{Lefschetz fibration} over $\mathbb{D}^2$ is a map from an exact symplectic manifold with corners $\pi: E \to \mathbb{D}^2$ which is $I_E$-holomorphic and satisfies some assumptions on behaviour near the boundary and the structure of critical points, as follows:

\begin{enumerate}
\item Transversality to $\partial \mathbb{D}^2$.

 At every point $x \in E$ such that $y= \pi(x) \in \partial \mathbb{D}^2$, we have $T  \mathbb{D}_y^2 = T \partial \mathbb{D}_y^2+D\pi (T E_x)$.

This implies that $\pi^{-1} (\partial \mathbb{D}^2)$ is a boundary stratum of $E$ of codimension 1, and we call it the vertical boundary of $E$, denoted by $\partial^v E$. The union of boundary faces of $E$ not contained in $\partial^v E$ is the horizontal boundary of $E$, denoted $\partial^h E$.

\item Regularity along $\partial^h E$.

 If $F$ is a boundary face of $E$ not contained in $\partial^v E$, then $\pi|_F:F \to \mathbb{D}^2$ is a smooth fibration. 

This implies that any fiber is smooth near its boundary.

\item Horizontality of $\partial^h E$ with respect to the symplectic connection. 

At any point $x$ of $E$, we have $TE_x^v= \operatorname{ker}(D \pi_x)$. Away from critical points, the fact that $\pi$ is $I_E$ holomorphic implies that the symplectic complement $TE_x^h$ of $TE_x^v$ is transverse to it (and so defines a connection).  We call $TE_x^h$ the \emph{horizontal tangent space} at $x$. We require that for all $x$ in any boundary face $F$ in $\partial^h E$ the horizontal $TE_x^h$ is contained in $TF_x$.

\item Lefschetz singularities. 

We require that the critical points of $\pi$ are generic (also called nondegenerate) and locally integrable. This means
that $I_E$ is integrable in a neighborhood of the set of critical points $Crit(\pi)$, and that $D \pi$ (seen as a section
of the bundle $Hom_{\mathbb{C}}(TE; \pi^{*}T\mathbb{D}^2)$ of complex linear maps) is transverse to the zero-
section. The second condition is equivalent to saying that the complex Hessian $D^2 \pi$
at every critical point is nondegenerate. In addition, we will assume that there is at
most one critical point in each fiber, so that, denoting the set of critical values of $\pi$ by $Critv(\pi)=\pi (Crit(\pi) )$, we get that the projection $Crit(\pi) \to Critv(\pi)$ is
bijective; this last assumption is for convenience only, and could easily be removed.

\end{enumerate}

\end{dfn}

 Nondegeneracy of critical points implies that they are isolated, so $Crit(\pi)$ is a finite
subset of $\operatorname{int}(E)$, and similarly $Critv(\pi)$ a finite subset of $\operatorname{int} (\mathbb{D}^2)$. Locally near each
critical point and its value, one has holomorphic coordinates in which $\pi$ becomes
the standard quadratic map
$Q(x) = x_0^2+ \ldots + x_{n}^2$. 
Generally $\omega_E$ will not be standard in these coordinates. However, one can find a
deformation of the fibration which is well-behaved along $\partial_h E$ (and which in fact
is local near the critical point), such that at the other end of the deformation the
K\"ahler form becomes the standard form in a given holomorphic Morse chart. This
increases the importance of the following basic model (see also Figure \ref{LefSing}):
 
 \begin{ex}\label{toy} (\cite[15.4]{Book}). Let $Q :\mathbb{C}^{n+1} \to \mathbb{C}$ be the quadratic function $Q(z_0, \ldots, z_n)=z_0^2+\ldots+z_n^2$, and $k :\mathbb{C}^{n+1} \to \mathbb{R}_{\geq 0}$ the
function $k(z) =\frac{ (|z|^4 -|Q(z)|^2)}{4}$. 
For some fixed $r, s > 0$ define
$ E = \{z \in \mathbb{C}^{n+1} : |Q(z)| \leq   r;  k(x) \leq s \}$, 
and equip it with the restriction of the standard symplectic form on $\mathbb{C}^{n+1}$, its standard primitive $\frac{i}{4} (z d \overline{z} - \overline{z} d z)$, the given complex structure $I_E = i$, and the map
$\pi : E \to r \mathbb{D}^2$ obtained by restricting $Q$. 
 
 The boundary faces are
$\partial_v E = \{z \in  E : |\pi(z)| = r\}$, $ \partial_h E = \{z \in  E : k(z) = s\}$.
The cutoff function $k$ is chosen so as to make $TE^h$ parallel to $\partial_hE$. To see that, one
notes that $TE^h_z$ is generated over $\mathbb{C}$ by $(\nabla Q)_z = 2 \overline{z}$, and checks that $dk_z(\overline{z}) = 0$,
$dk_z(i \overline{z})$ = 0. Each nonsingular fiber $E_w$, $w \neq 0$, is symplectically isomorphic to
the subset $B^*_s S \subset T^* S^n$ consisting of cotangent vectors of length (in the standard
metric) at most $\sqrt s$. Explicitly, $B^*_s S^n = \{(p, q) \in \mathbb{R}^{n+1} \times S^n : \langle p, q \rangle = 0, |p|^2 \leq s \}$, with the
symplectic form $dp \wedge dq$, and an isomorphism $E_w \to  B_s S^n$ for $w > 0$ is given by
$\phi_w(z) = (- \operatorname{Im}(z) |\operatorname{Re}(z)|,  \operatorname{Re}(z)|\operatorname{Re}(z)|^{-1})$. 
We emphasize that under this isomorphism $\phi_w(z)$ the function $k(z)$ becomes $|p|^2$, a fact we will use when we return to the geometry of this example in sections \ref{sectVC} and \ref{sectConstr}.

Unfortunately, while the  Liouville (radial) vector field does point
outward along $\partial E$, it is not true that $E$ is weakly $I_E$-convex (the fibers are, but
not the total space). Hence this is not quite an example of an exact Lefschetz
fibration as defined here, even though from a purely symplectic viewpoint, it has all
the desired features. One could change $I_E$ to remedy this, as we will do in Section \ref{sectConstr} when faced with a similar issue while constructing a Lefschetz fibration on $T^*S^{n+1}$.  We do not pursue this here, since we will use this example $E$  only as a local
model.
\end{ex}

The existence of this local normal form is a
consequence of the holomorphic Morse Lemma (more precisely, the statement is
that for any choice of holomorphic coordinates on the base, one can find coordinates
on the total space in which  $\pi= Q$). The deformation which allows one to make the
symplectic structure standard in such coordinates is constructed in \cite{S2}, Lemma
1.6. 

The total space $E$ of a Lefschetz fibration is an exact symplectic manifold with corners. As mentioned in the end of last section, after smoothing the corners and completing, one obtains a Liouville manifold which we will call $\hat{E}$. We will occasionally say that a Lioville manifold $G$ admits a Lefshcetz  fibration if $G$ is Liouville isomorphic to $\hat{E}$ for some Lefschetz fibration $\pi$ with total space $E$.  

\end{section}

\begin{section}{Vanishing paths and cycles, Dehn twists, Lefschetz thimbles and matching cycles.}\label{sectVC}

\begin{subsection}{Vanishing paths and vanishing cycles.}

\begin{wrapfigure}{r}{0.43\textwidth}
\centering
\includegraphics[width=0.43\textwidth]{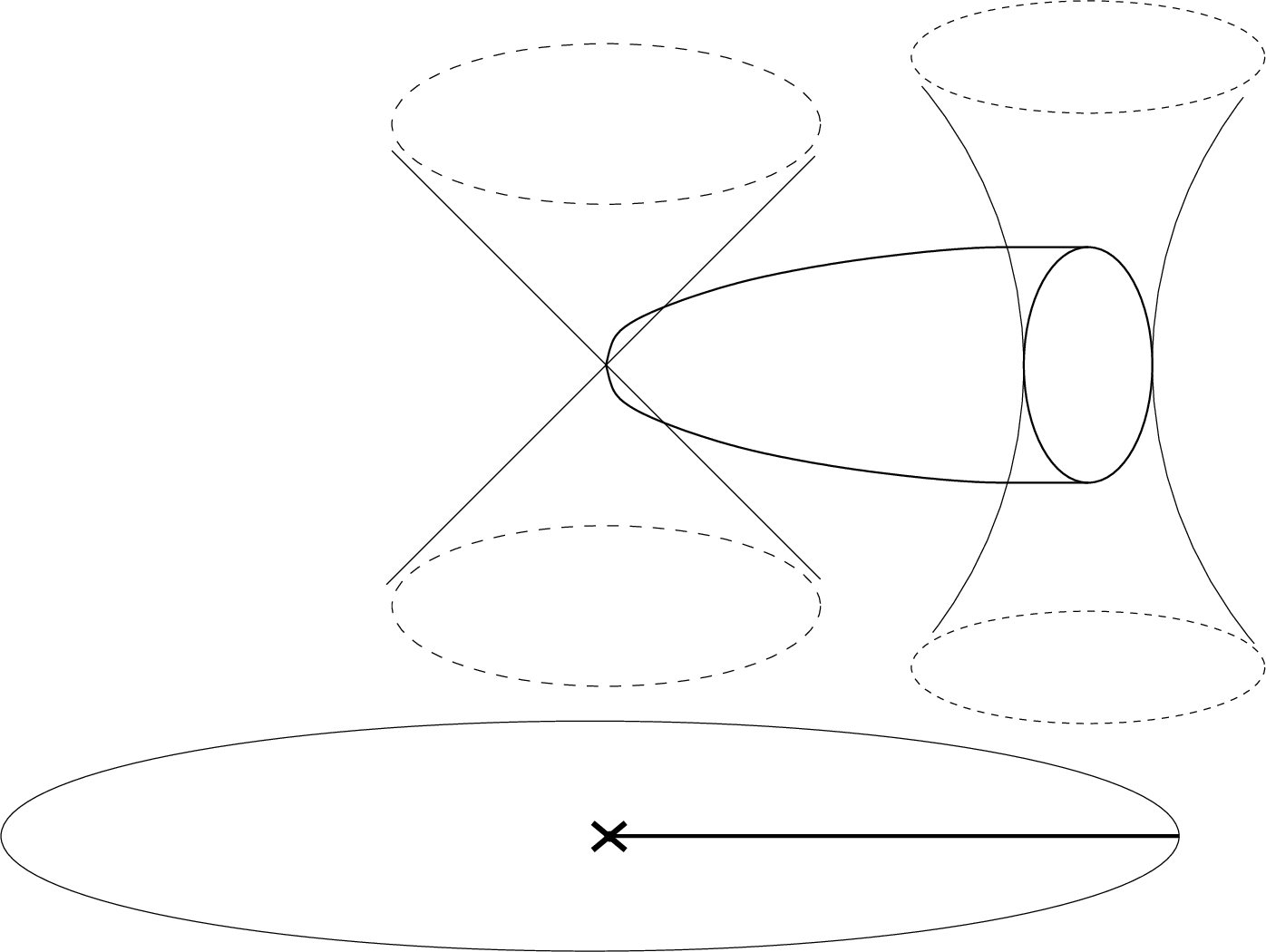}
\caption{Model Lefschetz singularity.}
\label{LefSing}
\end{wrapfigure}

For an (exact) Lefschetz fibration, we call an embedded curve  $\gamma:[0,1] \to \mathbb{D}^2$  a \emph{vanishing path}  if  it avoids critical points except at the end, i.e.  $|\gamma^{-1}(Critv (\pi))|={1}$. To each such path we can associate its \emph{Lefschetz thimble} which is the unique embedded Lagrangian
$(n + 1)$-ball in $E$ satisfying $\pi (\Delta_{\gamma})= \gamma([0,1])$. The boundary
$V_{\gamma} = \partial \Delta_{\gamma}$, which is a Lagrangian sphere in  $E_{\gamma(0)}$, is called the \emph{vanishing cycle} of $\gamma$. Since it bounds a Lagrangian disc in $E$, any vanishing cycle is automatically exact. We refer to section 16 of \cite{Book} for the proof that such a thimble exists and is unique. What is relevant for us is Remark 16.4,  which states that  the vanishing cycle of any piece $\gamma |[t_0; 1]$ is related to that of the whole by parallel transport. The precise statement is as follows. Remember that away from $Crit (\pi)$, we have a connection defined by the symplectic orthogonals to the tangent space of the fiber of $\pi$. We denote the parallel transport map of this connection along a path $\alpha$  in $\mathbb{D}^2 	\setminus Critv(\pi)$ by $h_\alpha$. Then $V_{\gamma}=h^{-1}_{\gamma | [0, t_0]} (V_{\gamma | [t_0, 1]})$. 
  Moreover, the vanishing cycle comes with an isotopy class of  \emph{framings}  - diffeomorphism $v: S \to V$ of the standard sphere $S$ (near the critical point this comes from the diffeomorphism with the sphere in tangent space at the critical point, and is then moved to the vanishing cycle itself by the parallel transport  diffeomorphism along the vanishing path). This isotopy class is part of the data of the vanishing cycle.


\begin{ex}\label{lef}
Take the model $\pi: E \to \mathbb{D}^2$ as defined in Example \ref{toy}. As mentioned before, this
is not quite a Lefschetz fibration, but the missing condition (lack of holomorphic
convexity) is irrelevant for the present purpose. The only critical value is 0, and
for any vanishing path $\gamma$, the Lefschetz thimble can be explicitly determined: $\Delta_\gamma = \bigcup \limits_{0\geq t \geq 1}\sqrt{\gamma(t)} S^n$.

Here 
$\sqrt z S^n=  \{ x \in  C^{n+1}: x = \pm \sqrt z y$  for some $y \in  S^ n \subset R^{n+1}\}$.
To see that this is the case,  one uses the function $k$ from Example \ref{toy}, which is
unchanged under parallel transport, and observes that $k^{-1}(0)$ is precisely the union
of the subsets $\sqrt z S^n$ for all $z \in \mathbb{D}^2$. In the identifications of the fibers of the model fibration in Example \ref{toy} with sphere cotangent bundles, these spheres $\sqrt z S^n$  are the zero sections. These spheres are what we will refer to as the \emph{belt spheres} in what follows.

\end{ex}

\end{subsection}

\begin{subsection}{Dehn twist.}

The existence of symplectic connection also allows us to associate the \emph{monodromy map} to any closed loop $\beta$ in the base avoiding all critical values. Namely, parallel transport defines a map $\mu_\beta$ from the fiber over $\beta(0)$ to itself; this map is symplectomorphism and changing the loop $\beta$ by a homotopy (avoiding the critical values of $\pi$) changes $\mu_\beta$ by an isotopy.

 Moreover, we have the folowing description of the monodromy map.

Suppose we have a vanishing path $\gamma$ connecting a reference point $\gamma(0)$ to a critical value $\gamma(1)$. We can create a loop $\beta$ from $\gamma(0)$ to itself, by going along $\gamma$ untill we are in a small contractible neighborhood of $\gamma(1)$, then going once counterclockwise  around $\gamma(1)$ while staying in that small neighborhood, and then going back along $\gamma$ traversed in reverse. The homotopy class of the resulting loop $\beta$ (for homotopies relative the basepoint, avoiding $Critv(\pi)$), depends only on the homotopy class of $\gamma$  (relative endpoints, and avoiding $Critv (\pi)$, see \cite{Book}[Section 16c]. It turns out that the isotopy class of the resulting monodromy symplectomorphism can be determined explicitely - it is represented by a \emph{Dehn twist} around the vanishing cycle  $V_\gamma$, an explicitly given homeomorphism supported near that $V_\gamma$. It is described in \cite{Book}[Section 16c] as well, and is named so because, if the real dimension of the fiber is 2, it coincides with the more familiar topological Dehn twist automorphism. In general, for any framed Lagrangian sphere $V$ in any symplectic manifold $M$ one can define the Dehn twist $\tau_V$ around that spehere - a symplectic automorphism of $M$ defined up to an isotopy.  We remark that we will not use the precise nature of the Dehn twist, but only it's effect on Floer cohomology (see Section \ref{sectDist}). 
\end{subsection}
 
\begin{subsection}{Matching cycles.}\label{secMatching}

Given a Lefschetz fibration $\pi$, we would like to collect the information contained in various vanishing paths and corresponding vanishing cycles. To begin, we want to have a collection of paths, one per Lefschetz singularity. This motivates the following definition.

\begin{dfn}\cite[Section 16d]{Book} Suppose the critical values of $\pi$ are $p_1, \ldots p_m$. Pick a base point $p$ on the boundary of $D$. A \emph{distinguished basis of vanishing paths} for $\pi$ is a set of $m$ vanishing paths $\gamma_i$, starting at $p$, that is $\gamma_i(0)=p$, and otherwise disjoint, and running into corresponding critical value, $\gamma_i(1)=p_i$. In addition, if one looks at the outgoing directions of $\gamma$s, that is  $(\gamma'_1(0),  \gamma'_2(0),   \ldots,  \gamma'_m(0))$, these directions should be arranged \emph{clockwise}.
\end{dfn}

Conversely, given a distinguished basis of vanishing paths in $D$ together with a collection of corresponding vanishing cycles one can reconstruct a Lefschetz fibration $\pi$.

\begin{lem}{\cite[Lemma 16.9]{Book}}\label{constrlem}  Fix a collection $(V_1, \ldots, V_m)$ of framed Lagrangian spheres in an exact symplectic manifold with corners $M$, a base point $p$ on the disc $D$ and a distinguished basis of vanishing paths $\Gamma=(\gamma_1, \ldots, \gamma_m)$. Then there exists a Lefschetz fibration $\pi: E \mapsto D$ with critical values endpoints of  paths in $\Gamma$ and an identification of the fiber of $\pi$ over $p$  with $M$ under which the (framed) vanishing cycles $V_{\gamma_i}$ correspond to $V_i$.  
\end{lem}

Moreover, these data of a distinguished basis of vanishing paths and corresponding vanishing cycles uniquely determine the fibration up to a deformation.

We now describe a \emph{matching cycle} construction. For an (exact) Lefschetz fibration $\pi$  consider an embedded path $\mu: [-1; 1] \to \operatorname{int} (\mathbb{D}^2)$ such that $\mu^{-1}(Critv(\pi)) = \{-1; 1\}$. We can split this into a pair of vanishing paths with the same starting
point, $\gamma _{\pm}(t) = \mu(\pm t)$ for $t \in [0; 1]$, hence get a pair of vanishing cycles $V_{\gamma \pm} \subset M =
E_{\mu(0)}$.  When these
two are equal (which is not going to be true on the nose in general, but suffices for the present applications)  $\Sigma_\mu = \Delta_{\gamma +} \cup  \Delta_{\gamma -}$
is a smooth Lagrangian submanifold of the total space $E$ (by definition of the Lefschetz
thimble, parallel transport along $\mu$ maps the intersections $\Sigma_\mu \cap \mu^{-1} (t)$ to each
other for all $-1 < t < 1$, which gives a local chart 
$(-1; 1) \times V_{\gamma \pm}$ around the overlap
$ \Delta_{\gamma +} \cap  \Delta_{\gamma -}= V_{\gamma \pm}$).
Being the result of gluing two balls along their boundaries, $\Sigma_\mu$ is necessarily a homotopy sphere. 
In the case when the framings of the $V_{\gamma{\pm}}$
are isotopic, it is a standard sphere differentiably. In fact, given a choice
of isotopy between the two framings, one can obtain a framing $\Sigma_\mu$. We will refer to $\Sigma_\mu$ as the \emph{matching
cycle} (see also Figure \ref{TSn1}).

We will apply this construction in the following context. Given a Lefschetz fibration $\rho: F \to \mathbb{D}^2$, we can choose some matching paths $\mu_i$. In the circumstances when the matching cycle construction above goes through, we obtain framed Lagrangian spheres $L_i$ in $F$. We then use the $L_i$'s to construct another Lefschetz fibration $\pi: E \to \mathbb{D}^2$ with vanishing cycles $L_i$. 

We note that in this case $E$ is an instance of a \emph{bifibration}. Bifibrations are discussed in some detail in section 15 of \cite{Book}, but we will not use their theory in any systematic way. 

\begin{rem}\label{remNoBdry} Observe that the Lefschetz thimble, vanishing cycle and matching cycle constructions do not rely on the almost complex structure and hence on weak convexity condition for the Lefschetz fibration. Moreover, for the purposes of these constructions, as long as the parallel transport involved in the cunstruction does not move the relevant vanishing cycle near the boundary of the  fibration, the horizontality condition (3) can be  omitted as well. We will use this to be able to talk about matching cycles in the statement of Proposition \ref{prop-SphereLF}.

\end{rem}

Finally we should note that the present discussion is somewhat simplified. Among other things,  one can define Lefschetz fibration over any Riemann surface with boundary, and give a more robust definition of matching cycles. Both of these and more can be found in the main reference for section 3 - Seidel's book \cite{Book}.

\end{subsection}

\end{section}

\begin{section}{The construction.} \label{sectConstrAll}

\begin{wrapfigure}[12]{r}{0.48\textwidth}
\centering
\includegraphics[width=0.48\textwidth]{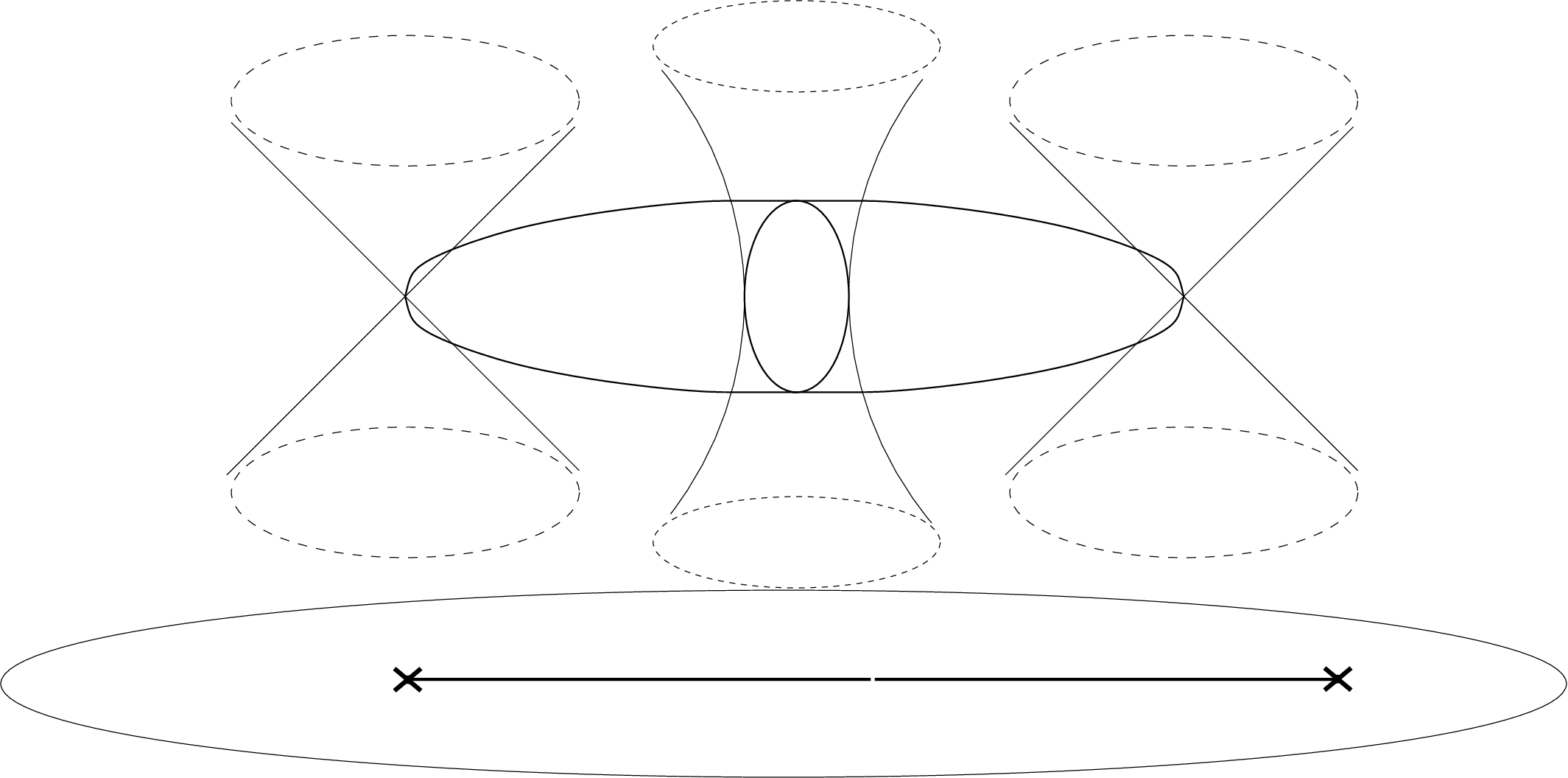}
\caption{Lefschetz fibration for $T^* S^{n+1}$.}

\label{TSn1}
\end{wrapfigure}

Our goal in this section is to construct two Liouville domains $W_0$ and $W_1$ and prove that they are diffeomeorphic. We will then show that they are not exact symplectomorphic in Section \ref{sectDist}.

We will build $W_0$ and $W_1$ by starting with a description of $T^* S^{n+1}$ as a Lefschetz fibration and modifying it in stages. This is depicted schematically in Figure \ref{InOut}, and will be explained in detail below. To build $W_0$ and $W_1$, at each stage we modify either the vanishing cycles of the Lefschetz fibration or its fiber. Lemma \ref{constrlem} then tells us that we can build the the corresponding total spaces.  Both $W_0$ and $W_1$ will be obtained in this manner.

We note that to make a compact part $E$ of $T^* S^{n+1}$ into a bona fide Lefschetz fibration requires making technical arguments aimed at insuring good behavior near $\partial E$. We provide them for the sake of completeness, but since ultimately $T^* S^{n+1}$ serves mostly as a \emph{topological} prototype for $W_1$ and $W_2$, these technicalities are somewhat tangential to the main development of this paper.

\begin{subsection}{The Lefschetz fibration on  $T^* S^{n+1}$.}
 
The description of $T^* S^{n+1}$ that we will give is as follows (see Figure \ref{TSn1}).

\begin{prop} \label{prop-SphereLF}There is a Lefschetz fibration $\pi:E  \mapsto D$  over an ellipse $D$ in  $\mathbb{C}$, with critical values $+1$ and $-1$ and with smooth fibers exact symplectomorphic to a disc bundle of $T^* S^n$ and vanishing cycles for any vanishing path identified with the zero section $S^n$, and such that the result of rounding corners of $E$ and completing to a Lioville manifold is Liouville isomorphic to $T^* S^{n+1}$.  The reference fiber $F$ of $\pi$  is  the total space of  a map  $\rho$ which is holomorphic and has only isolated Lefshcetz  singularities. Map $\rho$  two singularities with critical values $+1$ and $-1$, and the vanishing cycles of the main fibration $\pi$ are matching cycles for the straight matching path between the critical values of $\rho$. 

\end{prop}

\begin{rem} The map $\rho$ does not satisfy the conditions near the boundary that would make it a Lefschetz fibration. Nonetheless, because it is the same map as $\pi$ in one dimension lower, in light of Remark \ref{remNoBdry}  we can still talk about the matching cycles for paths that stay away from the boundary of the image of $\rho$. We could make $\rho$ into a genuine Lefschetz fibration at the expense of restricting it to a subdomain and deforming the complex structure. Since these arguments are essentially a repetition of the ones used for $\pi$, and we have no use for their outcome, we omit them. 
\end{rem}

 \begin{proof}

We work with cotangent bundles of spheres $T^* S^{n+1}$. These have the standard embedding into $\mathbb{R}^{2(n+2)}=T^*\mathbb{R}^{n+2}$, via $T^* S^{n+1}=\{(q,p)|  |q|=1, q \cdot p =0\}$, the derivative of the standard embedding of $S^{n+1}$ into $\mathbb{R}^{n+2}$. 
 In this model the standard symplectic form on $T^* S^{n+1}$  is the restriction of the symplectic form $\omega=d p \wedge dq$ by naturality.  However, for us a different model is going to be more convenient. 

Namely, consider the conic $\hat{M}=\{\Sigma z_j^2=1\}$ in $\mathbb{C}^{n+2}$.  This is the fiber over $1$ in the basic model of Lefschetz fibration see Example \ref{toy} (see also \cite{Book}, Example 15.9, or Lemma 1.10 in \cite{S2}). In terms of $x=\operatorname{Re} z$ and $y=\operatorname{Im} z$ in $\mathbb{R}^{n+2}$ it is given by $|x|^2-|y|^2=1$,  $x   \cdot y =0$, and hence  $\Phi: (x,y) \mapsto (-\frac{x}{ |x|}, y |x| )$ is a diffeomorphism to the sphere cotangent bundle. 

\begin{lem} The map $\Phi$ is an exact symplectomorphism between the conic $\hat{M}$ and the standard cotangent bundle of the sphere $T^* S^{n+1}$.

\end{lem}

\begin{proof}
We compute $d q_i=-\Sigma_j (\delta_{ij}  |x|^{-1} + x_j (-x_i) |x|^{-3}) d x_j$, so that $\Sigma_i p_i dq_i$ $=  -\Sigma_i y_i dx_i  - \Sigma_{i,j} x_j y_j x_i |x|^{-2} d x_i $ $= -\Sigma_i y_i dx_i,$ where we used $\Sigma_j x_j y_j =0$ along $M$ in the last equality. The standard primitive of the symplectic form on $\mathbb{C}^{n+2}$ is $\frac{i}{4}(z d \overline{z} - \overline{z} dz)= \frac{1}{2} (x dy -y dx)$. The difference between the pullback of $p dq$ computed above and the restriction of this primitive is $x dy+y dx = d \langle x, y \rangle=0$ on the conic, and the lemma follows.

\end{proof}

We note that the inverse map $\Phi^{-1}$ is given by $(q,p) \mapsto (-aq, \frac{1}{a}p)$, where $a^2-\frac{1}{a^2}|p|^2=1$, so $|p|^2=a^4-a^2$, $a\geq 1$.

We want to use the conic model from now on, but technically it does not fit with our definitions - we want to work with a compact exact  symplectic manifold with boundary, but the conic (as well as the cotangent bundle itself) is non-compact. This is remedied by taking a bounded part  where $|p|^2< s$. Note that $|p|^2$ is precisely $k(z)$ in Example \ref{toy}. Since we are on the fiber $Q(z)=1$,  via the symplectomorphism $\Phi$ above this corresponds to $|z|<(1+4s)^{\frac{1}{4}}$. We denote $(1+4s)^{\frac{1}{4}}$ by $r$, and assume that $s$ is large, and hence so is $r$.  Then the resulting manifold $M$ has an outward pointing Liouville flow and is weakly convex with respect to the standard complex structure. Moreover its completion is indeed $\hat{M}$.

We now build a Lefschetz fibration for this conic model. Consider the projection to  the last coordinate $\pi: \hat{M} \to \mathbb{C}$ sending $(z_0,\ldots z_{n+1})$ to $z_{n+1}$. Perhaps the easiest way to understand it is to note that the fiber over $\lambda$ is given by $z_0^2+\ldots+z_{n+1}^2=1-\lambda^2$, and so it is in fact a pull back via $\lambda \to 1-\lambda^2$ of the canonical Lefschetz local model fibration of Example \ref{toy} in one dimension lower. That is, there is a pull back diagram of holomorphic 
maps: 

\vspace {5mm}

$\xymatrix@1{ ( z_0, \ldots, z_{n} )\in \CC^{n+1} \ar[d]^Q & (z_0, \ldots, z_{n+1}) \in \hat{M}   \subset \CC^{n+2} \ar[l]_{Pr_{n+1}} \ar[d]^\pi  \\ z_0^2+\ldots + z_{n}^2=1-z_{n+1}^2 \in \CC  & z_{n+2} \in \CC  \ar[l]_(.3){1-\lambda^2} }$

\vspace {5mm}

In particular  $\pi$  has critical points when $z_0=\ldots=z_{n}=0$ and $1-z_{n+1}^2=0$ i.e. $z_{n+1}=1$ or $z_{n+1}=-1$; these critical points are nondegenerate; the corresponding critical fibers are conical.

Note that the maps in this diagram are not symplectic. However, when restricted to the fiber $F_{\lambda}$ of $\pi$, the map $Pr_{n+2}$ becomes symplectic, thus giving symplectomorphism of $F_\lambda$ with the fiber of $Q$ over $1-\lambda^2$ and the smooth fibers are exact symplectomorphic to the cotangent bundles of the sphere in one dimension lower, $T^* S^{n}$. 

As such they come with their zero sections - the belt spheres from Example \ref{lef}.  We  observe that for arbitrary  paths the parallel transport takes the belt spheres to each other. This follows from \cite [Lemma 6.12]{KS} which says that the union of belt spheres over any path is Lagrangian, and \cite[Lemma 16.3]{Book} which says that this only happens when the belt spheres are parallel transported to each other. Alternatively, we can note, similarly to \cite[5.1]{TDual} that  $\hat{M}$ is equipped with  $O(n+1)$ action (on the first $n+1$ coordinates in $\CC^{n+2}$) by symplectomorphisms, and an additional symplectic involution $(z_0, z_1, \ldots z_{n+1}) \mapsto (-z_0, z_1, \ldots z_{n+1})$. Both of these actions preserve fibers of $\pi$, and the belt sphere of each fiber is the unique orbit of the $O(n+1)$ action preserved by the involution (this is true even on the singular fiber where the belt sphere shrinks to a point). Since parallel transport is compatible with these actions, it follows that the belt spheres are taken to each other by it, and, in the case of vanishing paths,  the unions of belt spheres make up the corresponding Lefschetz thimble. Moreover, just like in the case of the model Lefschetz fibration, the corresponding vanishing cycle has standard framing. This means that for the Lefschetz fibration at hand the naive matching cycle construction from Section \ref{secMatching} works for any matching path.

 Note that if we take a reference fiber above $\lambda=0$ and straight line vanishing paths (the intervals $[-1, 0]$ and [0,1]),  the union of the corresponding two thimbles is the zero-section $S^{n+1}$ of the total space (see Figure \ref{TSn1}).

We now restrict to $M$, that is the part of $\hat{M}$ with $|z|<r$. The part of the fiber of $\pi$ over a given $z_{n+1}$ where $|z|<r$ is given by $|z_0|^2+\ldots +|z_n^2|^2 \leq r^2-|z_{n+1}|^2$, and on this fiber $z_0^2+\ldots +z_n^2 = 1-z_{n+1}^2$. From Example \ref{toy}, we get that this part of the fiber is symplectomorphic to the disc cotangent bundle of an $n$-sphere of radius squared  $w=\frac{1}{4}((r^2-|z_{n+1}|^2)^2- |1-z_{n+1}^2|^2)= \frac{1}{2} (\frac{r^4-1}{2}-  x_{n+1}^2 (r^2-1) - y_{n+1}^2 (r^2+1) )$. 
In particular, the piece $M$ of $\hat{M}$ where $|z|<r$ projects to the inside of the ellipse $w\geq 0$ with foci at $z_n=1$ and $z_n=-1$ and axes of lengths $\sqrt{\frac{r^2+1}{2}}$ and  $\sqrt{\frac{r^2-1}{2}}$.   We take a smaller ellipse  $D=\{ w \mid w \geq \epsilon^2 \}$ and restrict our fibration to it and to the subset of the total space which corresponds to disc bundles of radius $\epsilon$ over the corresponding fibers.  We take  $\epsilon$  small compared to $s$ and $r$ but big compared to $1$. We call the total space of the resulting fibration with singularities $E$, and keep calling the restriction of the fibration $\pi$. 

We want to see that $E$ is a Liouville subdomain of $\hat{M}$, that is that the Liouville flow of $\hat{M}$ is everywhere outward pointing along the boundary of $E$. This is obvious for the horizontal piece of the boundary, since the Liouville flow of $\hat{M}$ restricts on the fibers of $\pi$ to the Liouville flow on the disc cotangent bundle of $S^{n}$ and is transverse to the boundary of that  disc cotangent bundle. For the vertical part consider the singular fiber of  our main fibration $z_0^2+\ldots+z_{n+1}^2=0$; on it the Liouville flow is the radial flow from the origin, which is clearly transverse to the subset of $\mathbb{C}^{n+2}$ defined by  $w = \epsilon^2$.  If we take $r$ large enough, the Liouville vector field on our $\hat{M}$, the fiber above $1$, will be a small perturbation of the Liouville vector field on the singular fiber, and hence also transverse to the vertical boundary. Alternatively, one can compute the Liouville vector field explicitly. This computation shows that the required transversality holds for any $r$ (and not only large ones).  

The map $\pi: E \mapsto D$ is almost, but not quite, a Lefschetz fibration.  We claim that the only issue is that the boundary of $E$ is not convex for the complex structure $E$ inherits from $\hat{M} \subset \CC^{n+2}$. Ignoring this for a moment, we check conditions (1)-(4) from Definition \ref{def-LF}. All of them except horizontality of $\partial^h E$ are obvious. For this last one, we argue as in  Example \ref{lef}. Namely, each fiber $F_w$ of $\pi$  is symplectomorphic via $\psi_w$ to $T^*_\epsilon S^n$, the disc cotangent bundle of $S^n$ of size $\epsilon$. The action of $O(n+1)$ restricts to $E \subset \hat{M}$ and to $F_w$ and under $\psi_w$ the orbits of this action are subsets of $T^*_\epsilon S^n$ where the cotangent vectors are of constant length, and the boundary of $F_w$  is the orbit where theis length is $\epsilon$.  The parallel transport takes orbits to orbits, it takes the belt spheres to belt spheres and is a symplectomorphism, hence it must take the boundaries of various $F_w$ to each other, as desired.

Finally, to make $E$ into a genuine Lefschetz fibration, we modify the complex structure on it.

Denote the complex structure inherited from $\CC^{n+2}$ by $I$, and the standard complex structure $D$ inherits from $\CC$ by $I_{std}$. We construct an almost complex structure $I_E$ on $E$ which makes the boundary $\partial E$ weakly convex, and keeps the projection $\pi$ holomorphic. 

Observe that the fact that projection $\pi$ to $\CC$ is holomorphic for $I$ implies that the vertical part of the boundary $\partial^v E$ is already weakly convex for $I$. We shall modify $I$ near $\partial^h E$ to make $\partial^h E$ weakly convex as well. Namely, we pick a reference fiber $F$ of $E$ with its complex structure $I_F$. Note that $F$ is weakly convex for $I_F$. We pick a small open neighborhood $V$ and a small closed neighbourhood $N \subset V$ of $\partial F$ and use parallel transport along straight paths to map $V$ to $V_w \subset F_w$ $N$ and $N_w \subset F_w$ for all $w$ in the base ellipse $D$ of $\pi$. This provides a trivialization $\iota:  V \times D \mapsto E$ and we define an $J=\iota^* (I_F \times I_{std}$ on the image of $N \times D$. Finally, we construct $I_E$ by interpolating between $I$ on $\iota N \times D$ and $I$ on $E \setminus \iota( V \times D)$. 
With this new $I_E$ the map $E \mapsto D$ is a bona fide Lefschetz fibration.

We have therefore constructed an exact Lefschetz fibration which after rounding corners and completing gives Lioville manifold $\hat{M}$ which is Lioville isomorphic to $T^*S^{n+1}$, as wanted.

\end{proof}

 \end{subsection}
 
 \begin{subsection}{Constructing $W_0$ and $W_1$}\label{sectConstr}

  We now proceed to build Liouville manifolds $W_0$ and $W_1$ by modifying the above fibration on $T^*S^{n+1}$ in stages and applying Lemma \ref{constrlem} repeatedly.

\begin{figure}
\centering
\includegraphics[width=\textwidth]{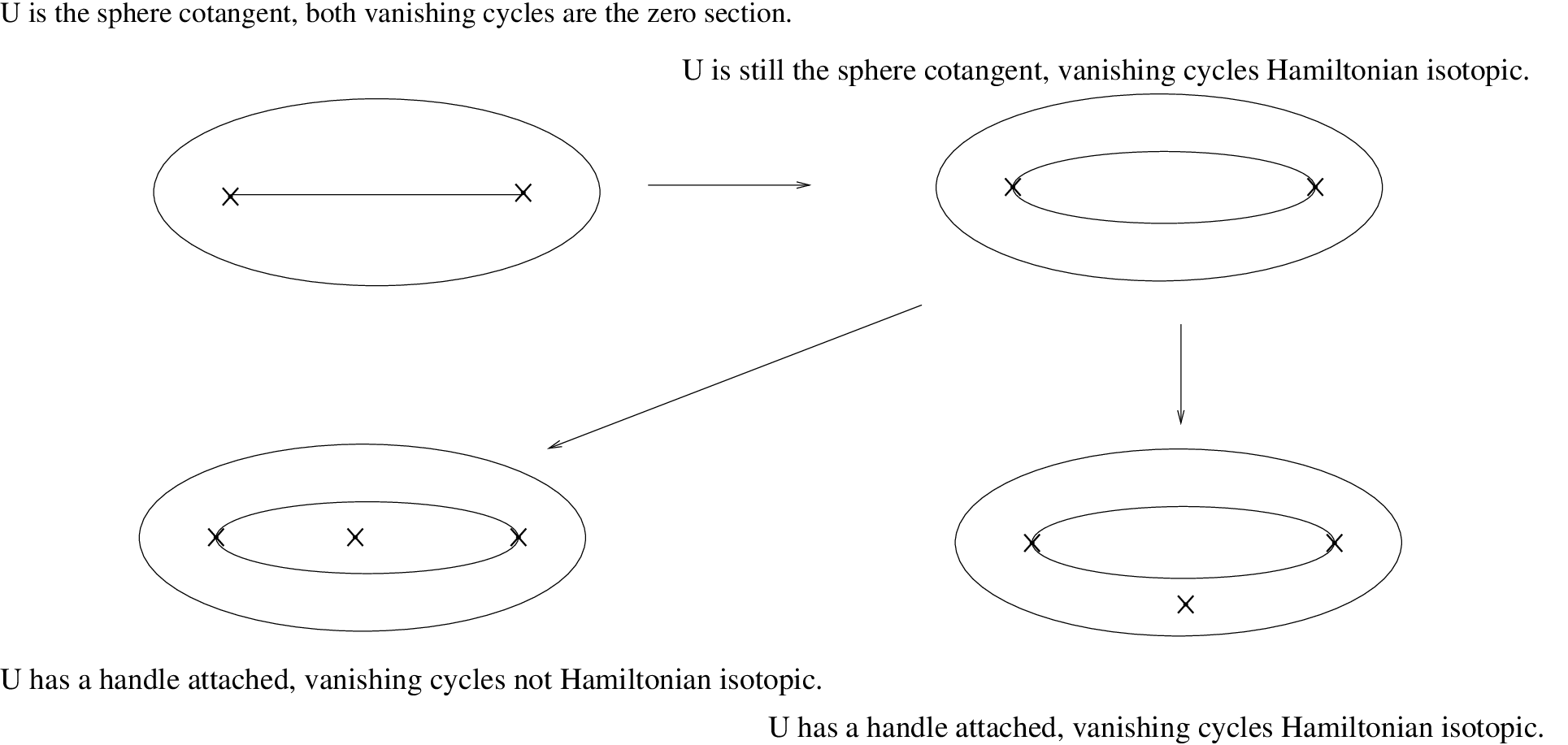}
\caption{Building $W_0$ and $W_1$.}

\label{InOut}
\end{figure}

This is schematically depicted in Figure \ref{InOut}. The first step is to change the matching paths for the auxiliary fibration $\rho$. The exact choice of paths is immaterial, for definiteness we can take the upper and lower semicircular arcs of $|z_{n+1}|=1$, which we will denote by $\alpha$ and $\beta$.  This produces matching cycles $A$ and $B$ in $T^*S^n$ that are framed Lagrangian isotopic to the ones we had before (the zero section).  If we were to build the total space $W$ still using $T^*S^n$ as fiber, but now using $A$ and $B$ as vanishing cycles, then by  the uniqueness counterpart of Lemma 16.9 in \cite{Book} after corner smoothing and completion, $W$ would be Liouville isomorphic to $T^* S^{n+1}$.

The second stage is more substantial. We change the fiber of the main fibration, i.e. the total space of the auxiliary fibration. The auxiliary fibration has two critical points $1$ and $-1$,  and for straight vanishing paths to a reference fiber at, say $-3i$ the vanishing cycles are the zero-sections of the fiber $T^* S^{n-1}$ over $-3i$.  We modify this by adding a third critical value at $-2i$ for $W_0$ and at $0$ for $W_1$, such that the vanishing cycle for a straight line vanishing path is equal to the belt sphere (in both cases).  Note that in the construction for $W_1$ the new critical value lies inside the disc encircled by the matching paths that define the Lagrangian spheres $A$ and $B$, while in the construction for $W_0$ it lies outside.  Applying Lemma \ref{constrlem}, we first get corresponding total spaces for the auxiliary fibrations, $U_0$ and $U_1$. These are symplectomorphic (topologically they are $T^* S^n$ with an $n$ handle),  but they contain different  matching cycles (see the left and the right diagrams in the bottom part of Figure \ref{InOut}, respectively). These matching cycles in $U_0$ and $U_1$ become vanishing cycles in the second pair of uses of Lemma \ref{constrlem}, which output the manifolds $W_0$ and $W_1$, correspondingly.

\end{subsection}
\begin{subsection}{$W_0$ and $W_1$ are diffeomorphic.}

We now show that $W_0$ and $W_1$ are the same as smooth manifolds. Later, in Section \ref{sectDist}, we will show that $W_0$ and $W_1$ are not exact symplectomorphic.

\begin{prop}\label{Same}

The smooth manifolds $W_0$ and $W_1$ are diffeomorphic.

\end{prop}

\begin{proof}

Consider the construction of $W_i$'s in more detail. Lemma 16.9 in \cite{Book} (which we have quoted  as Lemma \ref{constrlem}) describes the process of building the total space of the Lefschetz fibration as a sequence of surgeries - one first thickens the fiber by taking its product with $\mathbb{D}^2$ and then performs a series of handle attachments along the spheres given by the vanishing cycles (see also \cite{GompfStip}).  It is therefore sufficient to show that the matching cycles $A$ and $B$ are smoothly isotopic.

The matching cycle $A$ is the union of belt spheres of fibers over the matching path $\alpha$ and the matching cycle $B$ is the union of belt spheres of fibers over the matching path $\beta$ (plus the two critical points).  

The paths $\alpha$ and $\beta$ are isotopic in $\mathbb{D}^2$ (relative the endpoints). For $\alpha, \beta: [-1,1] \to \mathbb{D}^2$ we have the isotopy $\gamma:[-1,1] \times [0,1] \to \mathbb{D}^2$. We can take $\gamma$ to be symmetric, that is to satisfy $Re(\gamma(t, \cdot) = -Re (\gamma(-t, \cdot))$, and to stay inside the disc encircled by $\alpha$ and $\beta$.  Again, the exact choice of $\gamma$ is immaterial, but for definiteness we can take each $\gamma(\cdot, s)$ to be a uniformly parametrized circular arc through the two critical points and $i(1-2s)$.  For each time $s \in [0,1]$, if  the path $\gamma(\cdot, s)$ misses the  third critical point at the origin then it defines a matching cycle $\Gamma_s$. The only problem occurs when $\gamma_s$ hits the origin and the belt sphere over $\gamma(0,s)$ shrinks to a point (for our choice of $\gamma$ this happens at $s=\frac{1}{2}$). To remedy this, we push the $\Gamma_s$ off the zero section, see Figure \ref{Lifting}.

\begin{figure}
\centering
\includegraphics[width=\textwidth]{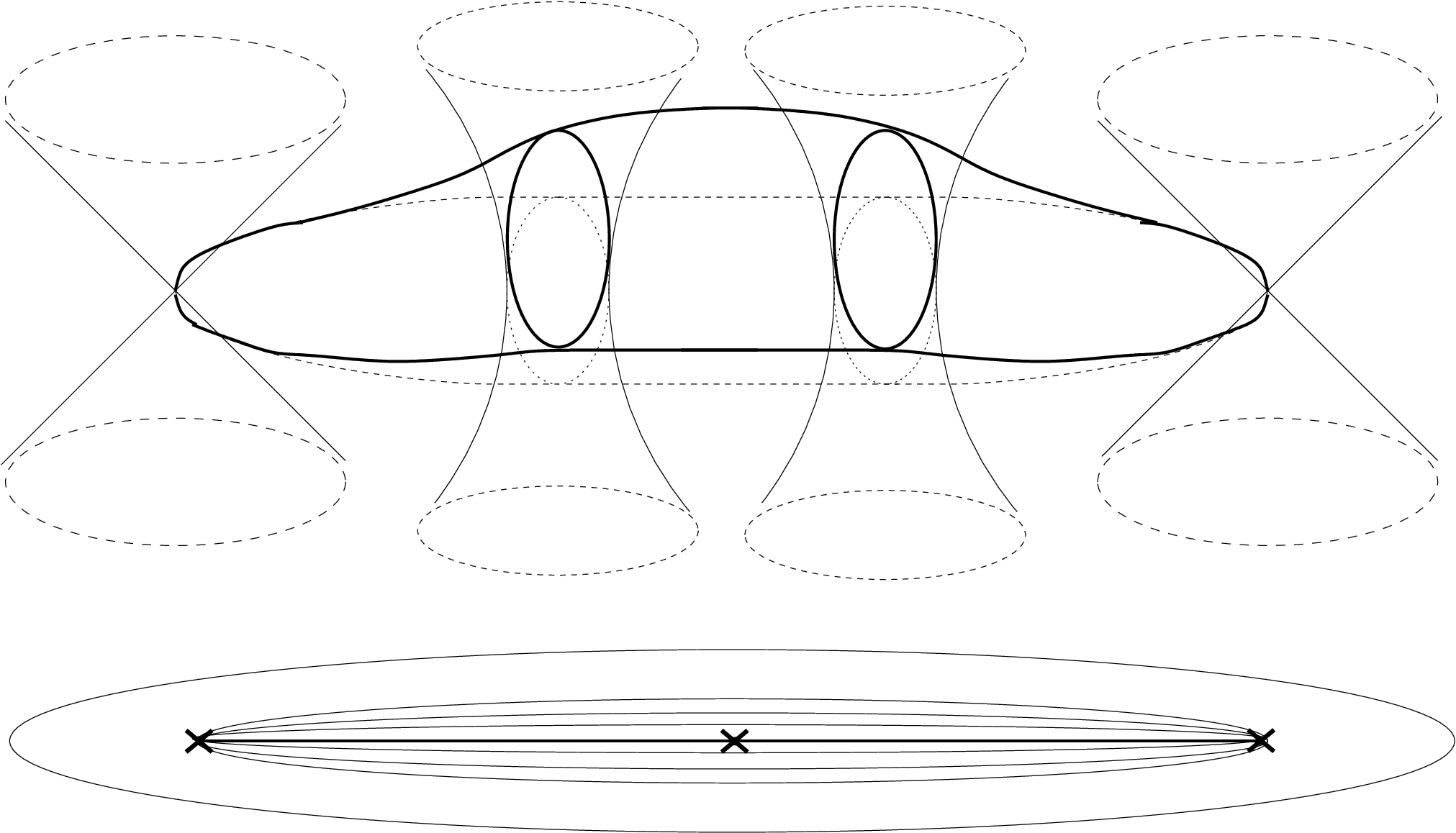}
\caption{Pushing the isotopy.}

\label{Lifting}
\end{figure}

The details are as follows. It is well-known that for even $n$ the sphere $S^{n-1}$ has a smooth  vector field $\delta$ with $|\delta|=1$ on it (in the standard metric). Take a smooth ``horizontal'' cutoff  function $h(t):[-1,1] \to [0, \nu<\frac{1}{2}\epsilon]$, which is zero near the endpoints and equal to $\nu$ near $0$; take also a smooth ``vertical'' bump function $v(s):[0,1]
\to [0,1]$.

Let $B$ be the disc enclosed by $\alpha$ and $\beta$ with small (contractible) neighborhoods of the critical points taken out.  Before we added the third critical point at zero, we had an auxiliary fibration $\rho:T^*S^n \mapsto \mathbb{D}^2$ with fibers $T^*S^{n-1}$. This was (smoothly) trivial over $B$.  Let  $$\Psi: B \times \{(v \in T^* S^{n-1} | |v|< \epsilon )\} \to \rho^{-1}(B)$$  be a trivialization taking belt spheres to belt spheres. 

 Define  for all $t$ such that $\gamma_s$ is in $B$, the set  $L(t,s)= \{ h(t) v(s) \delta(p) | p \in S^{n-1} \}$. This is a push-off of the zero section over $\gamma(t,s)$ in the direction of $\delta$.

We use the trivialization $\Psi$ to define $\hat{\Gamma}_s$ as the union of $L(t,s)$  over all such $t$ in the fiber  over $\gamma(s,t)$ together with the belt spheres of $\rho$ over the parts of $\gamma_s$ lying outside $B$ (these glue smoothly because $h(t)$ is zero near the endpoints). This is a sphere in the total space of $\rho$ - the push off of the whole matching cycle of $\gamma_s$ in the direction of $\delta$. Now the union of $\hat{\Gamma}_s$ over all $s$ gives an isotopy of the matching cycles $A$ and $B$ (recall that in the smooth category any non-ambient isotopy extends to an ambient one). Adding the third critical point at zero happens as a surgery on the belt sphere at zero, supported in its neighborhood. The spheres $\hat{\Gamma}_s$ stay away from the surgery region and hence persist in the manifold $U_1$, and hence can be used to define an isotopy between the matching cycles $A$ and $B$ in it.  If we take $\nu$ small enough, then the spheres $\hat{\Gamma}_s$ stay totally real, and hence the isotopy is an isotopy of framed spheres.

   This means that the total space  $W_1$ is the same smooth manifold as the space obtained by attaching handles to thickened $U_1$ along two copies of the matching cycle  $A$. But that is the same as the total space $W_0$.  This completes the proof.

 \end{proof}

  \begin{rem} In the case of $W_1$, if we denote by $L$ the matching cycle of the straight matching path from $-1$ to $0$ and by $R$ the matching cycle of the path from $0$ to $1$,   we see by Lemma 16.13 in \cite{Book} that the bottom matching cycle $B$ is obtained by a  Dehn twist around $R$ of $L$, and the top matching cycle $A$ is obtained by the inverse Dehn twist around $R$ of $L$. So $B$ is obtained from $A$ by the square of a Dehn twist around $R$.   For spheres of dimension 2 and 6, the square of the model Dehn twist is smoothly isotopic to identity (which of course implies that $A$ is isotopic to $B$). The case $n=2$ is Lemma 6.3 in \cite{knot}, and  the case $n=6$ can be handled somewhat analogously by using an almost complex structure on $S^6$ (this is an unpublished result of Giroux). The case of other $n$ appears to be open.   
  
  \end{rem}

\end{subsection}

\end{section}

\begin{section}{Wrapped Floer cohomology and wrapped Donaldson-Fukaya category.}

Wrapped Floer cohomology is an adaptation of the usual Lagrangian Floer cohomology to the context of non-compact Lagrangians in Liouville domains. The main additional concern in this situation is what to do when the Lagrangians in question have intersections at the boundary.

 As described in \cite[Section 1]{categor},  there are several possible solutions.  We shall describe two of them following \cite[Section 2]{categor}.

 Let $(M^{2n}, \theta)$ be a Lioville domain.  Recall that the Liouville flow $\kappa$ gives the collar neighborhood $\kappa:(-\infty, 0]\times \partial M  \to M$ with $\kappa^*\theta=e^s \alpha$ and $\omega=d\theta$, and this is extended to positive $s$ in the completion $\hat{M}$.

Following Abouzaid and Seidel \cite{AS}, we sometimes  use the coordinate $r=e^s$ on the cone. 

\begin{rem} This choice of coordinates for the infinite cone has the following slightly unfortunate consequence. Viewing, as before,  $\mathbb{C}$ as the completion of $\mathbb{D}^2$ with its standard Liouville structure, and denoting by $r_{polar}$ the standard radial polar coordinate on $\mathbb{C}$ as before,  we get $\omega= r_{polar} dr_{polar} \wedge d \theta = d(\frac{1}{2}r_{polar}^2 d \theta)$ so that $r=\frac{1}{2}r_{polar}^2$. In particular, functions which are linear in $r$ are  quadratic in $r_{polar}$. In discussing Hamiltonians on $\mathbb{D}^2$ and $\mathbb{C}$ we will work with the completion coordinate $r$, writing $r_{polar}$ whenever we refer to the standard polar radius. \end{rem}

\begin{subsection}{Floer cohomolgy for non-compact Lagrangians in Liouville domains.}

\begin{dfn} A Lagrangian $L \subset M$ is \emph{admissible} if it intersects $\partial M$ transversally and is exact, meaning $\theta|L$ is exact, $\theta|L=df$.   $L$  we require that $\theta|L$ vanishes on a neighborhood of $\partial L$. \end{dfn}

Near its boundary this makes $L$ a cone over the Legendrian submanifold  $\partial L$ of  the contact manifold $(\partial M, \theta|\partial M)$ and allows us to attach an infinite cone, extending $L$ to a Lagrangian $\hat{L} \subset \hat{M}$.

To define Floer cohomology $HF^*(\hat{L}_0,\hat{L}_1)$ we take a Hamiltonian $H \in C^\infty(M, \mathbb{R})$ which is everywhere positive and admits a smooth positive extension $\hat{H}$ to $\hat{M}$ such that $   \hat{H} (r,y)= r$ on the semi-infinite cone, and define $\hat{H}_l=l \hat{H}$ for any non-negative $l \in \mathbb{R}$.  Consequently, if $\hat{X}_l$ is the Hamiltonian vector field of $\hat{H}_l$, then on the cone $\hat{X}=(0, l R)$ where $R$ is the Reeb vector field of the contact one form $\theta|_{\partial M}$. 

We have the following lemma.

\begin{lemma} \cite[Lemma 2.1]{categor} \label{curvesel}
Let $\phi^t_R$ be the Reeb flow on $\partial M$. For any pair of Legendrians $(\partial L_0, \partial L_1)$, there
is an $\varepsilon > 0$ such that $\phi^t_R(\partial L_0) \cap L_1= \emptyset$ for all $t\in (0, \varepsilon)$.
\end{lemma}

Denoting by $\psi_{\hat{X}_l}$ the time one flow of $\hat{X}_l$,  consider the cochain complex $CF_*(\hat{L}_0, \hat{L}_1; \hat{H}_l)$ generated by the intersection of $\psi_{\hat{X}_l} (\hat{L}_0)$ and $\hat{L}_1$.   

If we take $l \in (0, \varepsilon)$, the intersection $\psi_{\hat{X}_l} (\hat{L}_0) \cap \hat{L}_1$ is compact by Lemma \ref{curvesel}. 

The Floer differential on this complex is given by counting holomorphic discs with boundary on $\psi_{\hat{X}_l}(\hat{L}_0)$ and $L_1$, or equivalently by counting solutions $u:[0,1] \times \mathbb{R} \mapsto \hat{M}$ of Floer equation $(du-X_l \otimes  dt)$ with boundary conditions on $L_0$ and $L_1$ (see \cite[Remark 10]{Ritter} ). If one chooses an almost complex structure $J$ of contact type (see Definition~\ref{contactJ}), then these solutions
stay away from the infinity by  \cite[ Section 15.22]{Ritter} (see also \cite[Section 3c]{Biased})  and we can define $HF(\hat{L}_0, \hat{L}_1)= H(CF_*(\hat{L}_0, \hat{L}_1; \hat{H}_l)$ in  the same way as for compact Lagrangians. As in the compact case, $HF(\hat{L}_0, \hat{L}_1)$ is independent of $H$ in this class and of $l \in (0, \varepsilon)$.

 Note that if $\partial L_0 \cap \partial L_1 =\emptyset$, we can actually take $l = 0$.   Finally, for $L_0 = L_1 = L$ we have 
 the usual reduction to Morse theory (this goes back to the classical computation of Floer (\cite{Floer2}[Theorem 3]; see also \cite{Book}[Sections 8c and 12e]), so that $HF^*(\hat{L},\hat{L}) =
H^*(\hat{L})$.

\end{subsection}

\begin{subsection}{Wrapped Floer cohomology of non-compact Lagrangians in Liouville domains.}

The wrapped Flooer cohomology is defined in a similar way, except instead of taking small $l$, we allow $l$ to grow and take a direct limit over the resulting system of cohomologies.  We note that wrapped Floer cohomology can be viewed as a Lagrangian (or open-string) version of symplectic cohomology introduced by Viterbo (\cite{V},\cite{V2}, see also \cite{Biased}).  It uses the same class of perturbations at infinity and the same direct limit procedure.

To define wrapped Floer cohomology we proceed as follows. Given two admissible Lagrangians $L_0$ and $L_1$, call $l \in \mathbb{R}$ admissible if there exists no Reeb chord in $\partial M$ of length $l$ starting at $L_0$ and ending at $L_1$.  Lemma~\ref{curvesel} showed that all $l\in (0, \varepsilon)$ are admissible, and a generic $l$ is admissible as well. 

 For an admissible $l$ the intersection $\psi_{\hat{X}_l}(\hat{L}_0)$ with $\hat{L}_1$ is contained in a compact subset of $\hat{M}$. For a generic $H$ the intersection $\psi_{\hat{X}_l}(\hat{L}_0)$ with $\hat{L}_1$ is in addition transverse. Denote by $CW_l (\hat{L}_0, \hat{L}_1)=CF(\hat{L}_0, \hat{L}_1; \hat{H}_l)$ the resulting Floer complex.

As in the previous section,  we can define $HW_l(\hat{L}_0, \hat{L}_1)= H(CF_l(\hat{L}_0, \hat{L}_1))$. Note that for small $l$, $HW_l({L}_0, {L}_1)=HF(\hat{L}_0, \hat{L}_1)$.

 Now, these cohomology groups come with continuation maps $HW_{l-}(L_0, L_1) \mapsto HW_{l+}(L_0, L_1)$ for  all admissible $l+ \geq l-$ (cf. \cite[Section 3.7]{Ritter}; this is similar to the case of symplectic cohomology discussed in  \cite{Biased}).

 Since direct limits are exact, we can define the wrapped Floer cohomology 
 as either the cohomology of the direct limit of the chain complexes or simply as the direct limit of the cohomology system:
  $$HW(L_0, L_1)= \lim_{l \to \infty} HW_l(L_0, L_1)$$
This is independent of the choices  made in constructing the Hamiltonians $\hat{H}_l$ and the choice of compatible $J$ of contact type. 

Given a triple of admissible Lagrangians $L_0$, $L_1$ and $L_2$, numbers  $l_{01}$, $l_{12}$, and $l_{02}=l_{01}+l_{12}$ with $l_{ij}$ admissible for the pair of Lagrangians $(L_i,L_j)$ and choosing the Hamiltonians so that $\hat{H}^{02}_{l_{02}}=\hat{H}^{01}_{l_{01}}+\hat{H}^{12}_{l_{12}}$, we get the product $CF(\hat{L}_1, L_2; \hat{H}^{12}_{l_{12}}) \otimes CF(\hat{L}_0, \hat{L}_1; \hat{H}^{01}_{l_{01}}) \mapsto CF(\hat{L}_0, \hat{L}_2; \hat{H}^{02}_{l_{02}})$ as in the unwrapped case by counting solutions of Floer equation on a disc with three boundary punctures and this passes through the directed system to give a product $HW(L_1, L_2) \otimes HW(L_0, L_1) \mapsto HW(L_0, L_2)$ (see and  \cite[Theorem 100]{Ritter} and \cite[Section 8a]{Biased}  for the case of symplectic cohomology).  All together, this makes admissible Lagrangians in $\hat{M}$ into objects of the wrapped Donaldsohn-Fukaya  category, with morphisms given by the wrapped Floer cohomology groups.

We should mention that in the work of Abouzaid and Seidel \cite{AS}  a slightly different version of wrapped Floer cohomology is given. There the authors are concerned with the $A_{\infty}$ structure on the chain complex, and need a finer model than the direct limit construction. We, on the other hand, are only concerned with the cohomology and product structure, so a cruder, simpler model is sufficient. The resulting cohomology is the same in both models (see Lemma 3.12 in \cite{AS}).

\end{subsection}

\begin{subsection}{Wrapped Floer cohomology in a Lefschetz fibration.} \label{sectWHF}

\medskip

\vspace{10mm}
  
    \parbox [t]{0.35\textwidth }{
\includegraphics [width =\linewidth ]{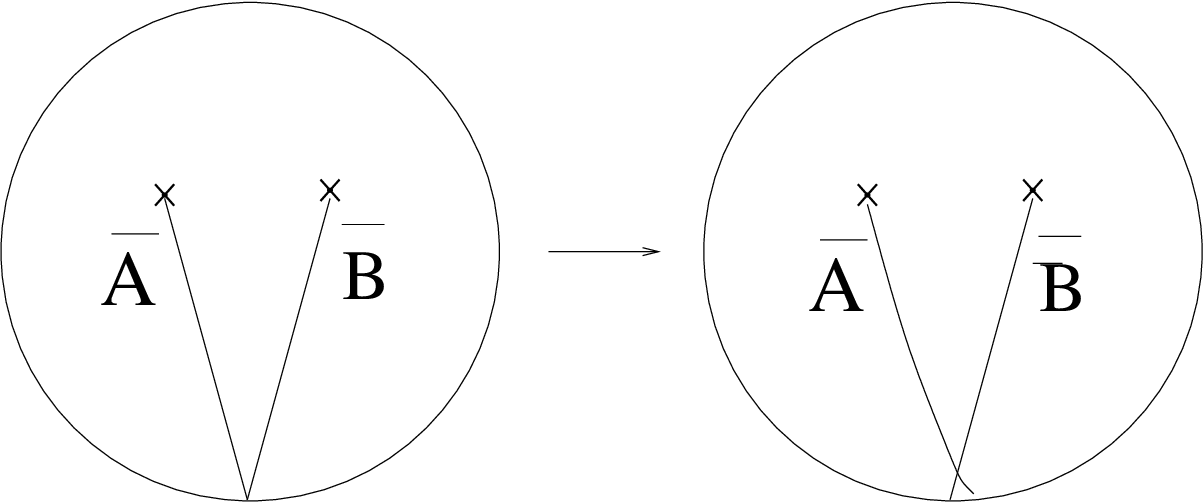}
\captionof {figure}{$HF(\overline{A}, \overline{B}) $}
}
\hfill
\parbox [t]{0.55\textwidth }{
\includegraphics [width =\linewidth]{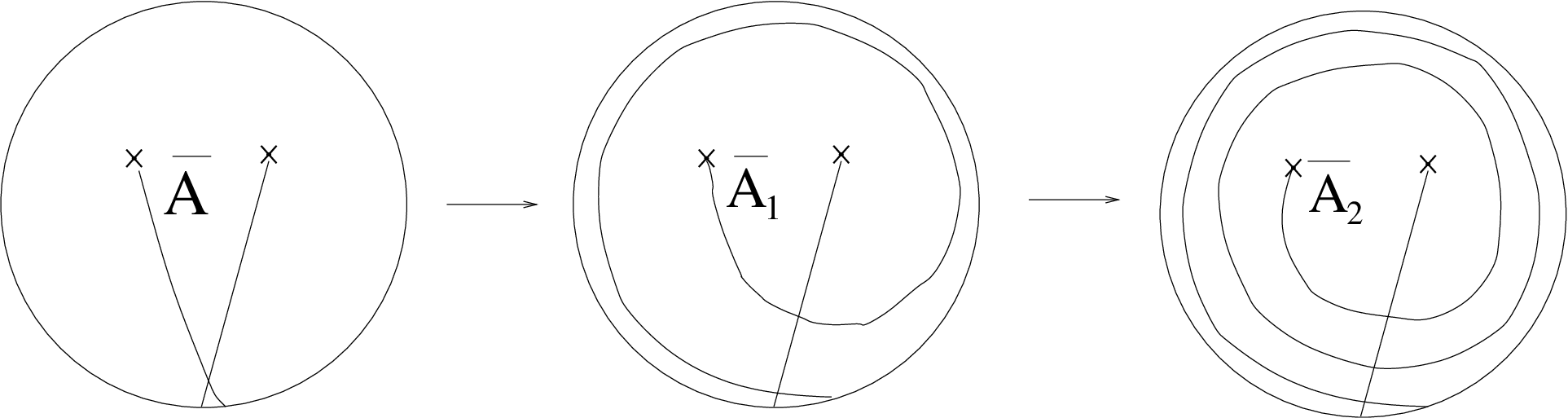}
\captionof {figure}{$HW(\overline{A}, \overline{B})=\lim HF_m(\overline{A}, \overline{B})$}
}

\vspace{5mm}

We now want to discuss an adaptation of these constructions to the case of Lefschetz fibration. This discussion is parallel to that of Section 3 in \cite{MS}.  We will consider wrapped Floer cohomology of thimbles $\Delta_\gamma$ for some vanishing path $\gamma$. The requirement that $\Delta_\gamma$ is admissible means that $\gamma$ is straight and radial outside of a compact set, that is for all $t$ large enough $\gamma(t)=|\gamma(t)| e^{i \alpha(\gamma)}$ for some constant angle $\alpha(\gamma)$. Then $\Delta_\gamma$ has a conical end modeled on the vanishing cycle $V_\gamma$, Legendrian in $\partial M$.

To make our construction of wrapped Floer cohomology adapted to the Lefschetz fibration, we modify it by considering a different family of  Hamiltonians.

Namely,  consider a Hamiltonian $H_b$ on the base disc of the fibration which is radial, zero in the interior of the disc and of slope one in $r$ near the boundary: $H_b(r e^{i \theta})=\eta(r)$, for $\eta(r):[0,1] \to [0, \infty)$ with $\eta(r)=0$ for $r<\sqrt{2}-2\epsilon$, and $\eta'(r)=1$ for $r>1-\epsilon$. Further, given a Lefschetz fibration $\pi : E\mapsto \mathbb{D}^2$,  define $H$ to be the pull back of $H_b$ under $\pi$. Then $H_b$ extends to all of $\mathbb{C}$ linearly in $r$, and $H$ extends linearly to $\hat{H}$ on the part of the infinite cone attached to the vertical boundary of the fibration.

The Hamiltonian $H_b$ generates a flow which is supported near the boundary of $\mathbb{D}^2$ and acts by rotations. We have the following observation:

\begin{obs} \label{lift}Let $X_H$ be the Hamiltonian vector field of $H$ generating the flow  $\phi$ . Let $X_{H_b}$ be the Hamiltonian vector field  of $H_b$ generating the flow $\phi_b$ . Then  $X_H$ propotional to the horizontal lift with respect to the symplectic connection of the fibration of $X_{H_b}$, and $\phi$ is a reparametrization of the lift of $\phi_b$. 
 \end{obs}

\begin{proof} For any  tangent vector $v \in TE_x$, we compute:  $H= \pi^* H_b$, $dH= \pi^* d H_b$,  $dH(v)= \pi^* d H_b (v)= d H_b(\pi_{*} v) $, so that
\begin{equation}
\nonumber  (*) \hspace{2cm} \omega(X_H, v)= \omega_{st} (X_{H_b}, \pi_*(v))
\end{equation}
  By taking vertical $v$, we see that $X_H$ is in the symplectic complement of the fiber. Further, taking $v=X_H$ we get $0=\omega(X_H, X_H)= \omega_{st} (X_{H_b}, \pi_*(X_H))$, so $\pi_*(X_H)= c X_{H_b}$ for some $c \in \mathbb{R}$.

 Hence $X_H$ is the  horizontal lift of $X_{H_b}$ with respect to the symplectic connection and the observation follows.

\end{proof}

To understand the proportionality constant $c$ we apply the same equation $(*)$ from the the proof of Observation \ref{lift} to $v= I_E(X_H)$. Denoting by $g(\cdot, \cdot)=\omega(\cdot, I_E(\cdot))$ and $g_{b}(\cdot, \cdot)=\omega(\cdot, I_{st}(\cdot))$ the metrics on $E$ and $D$ respectively, and using the fact that $\pi$ is holomorphic, we get  $g(X_H, X_H)= c g_b(X_{H_b}, X_{H_b})$. Hence $c$ is the ratio between the length of $X_{H_b}$ and the length of its horizontal lift. When constructing Lefschetz fibrations via Lemma \ref{constrlem}, we can ensure that $c$ is $1$ everywhere. Since we will apply wrapped Floer cohomology only to Lefschetz fibrations  constructed in this way, we may assume this is the case (more generally this can be achieved by a suitable deformation).

Correspondingly,  the result of flowing a Lefschetz thimble by $X_H$ is again a Lefschetz thimble,  namely $\phi(\Delta_\gamma)=\Delta_{\phi_b(\gamma)}$.  
The time 1 flow of $\hat{H}$ wraps the thimbles exactly once around $\partial\mathbb{D}^2$ at the boundary.  

We will use the following notation. We consider  vanishing paths $\alpha$ and $\beta$ with vanishing cycles $A$ and $B$ and thimbles $\overline{A}$, and  $\overline{B}$ correspondingly. We want to define Floer cohomology $HF(\overline{A}, \overline{B})$ and wrapped Floer cohomology $HW(\overline{A}, \overline{B})$ for the thimbles.

 We consider a family of Hamiltonians $\hat{H}^m= (m+\delta)\hat{H}$, $\hat{H}^m_b= (m+\delta)\hat{H_b}$
for positive integers $m$ and a fixed positive small real $\delta$. We denote the flow of $\hat{H}^m$ by $\phi^m$, and $\hat{H}^m_b$ by $\phi^m_b$. We choose $\delta$ so that $\phi^m_b(\alpha)$ does not intersect $\beta$ at the boundary (this is true for all $m$ if and only if it is true for $m=1$).

All the intersections of the thimbles $\overline{A}$ and $\overline{B}$ lie in the fibers over the intersection points of their paths in the interior of the disc, and are given by intersections of the corresponding vanishing cycles (for the case of the thimble intersecting with itself, this may require an additional compactly supported isotopy of the vanishing path to make it transverse to itself, which corresponds to a Hamiltonian isotopy of the thimble and does not affect the Floer cohomology).  After a small perturbation of the symplectic connection on the regular part of the fibration,  we can assume that the relevant  vanishing cycles intersect each other transversely. This perturbation is a compactly supported exact change of the symplectic form on $E$ and does not change the Floer cohomologies in the base or the fiber.  Then for a pair of thimbles $\overline{A}$ and $\overline{B}$
 we have for each $m$ the usual Floer-theoretic finitely generated chain complex  $CF_m(\overline{A}, \overline{B}) := CF(\psi_m(\overline{A}), \overline{B})$, with a differential obtained by counting pseudoholomorphic strips, which computes the Floer cohomology $HF_m(\overline{A}, \overline{B}):= HF(\psi_m(\overline{A}), \overline{B})$ in the usual way.  More precisely, the convexity conditions in the definition of exact symplectic manifold with corners ensure that such holomorphic curves stay in a compact part of the Lefschetz fibration, and hence their moduli spaces are compact, since bubbling is excluded by the exactness assumptions, and the rest of analysis is the same as in the compact case (just like in the case of the Liouville domains) and goes back to Floer \cite{Floer}. This is an extremely simplified (ungraded, with $\mathbb{Z}_2$ coefficients) cohomology-level version of the setup in \cite{Book} .

  Again as usual, this comes with a product structure $HF_{m_1}(\overline{A}, \overline{B})  \otimes HF_{m_2}(\overline{B}, \overline{C}) \to HF_{m_1+ m_2} (\overline{A}, \overline{C})$ obtained by counting holomorphic triangles. More precisely,  the composition of the two Hamiltonian isotopies $\phi^{m_1}$ and $\phi^{m_2}$ is a Hamiltonian isotopy $\phi^{m_1, m_2}$, which is generated by $(m_1+m_2+ 2\delta) \hat{H}$ and is close to $\phi^{m_1+m_2}$. We have the standard Floer product $HF_{m_1}(\overline{A}, \overline{B})  \otimes HF_{m_2}(\overline{B}, \overline{C}) \to HF (\phi^{m_1, m_2}(\overline{A}), \overline{C})$ and the small hamiltonian isotopy induced by $\delta \hat{H}$ induces canonical isomorphism $HF_{m_1+ m_2} (\overline{A}, \overline{C}) = HF (\phi^{m_1, m_2}(\overline{A}), \overline{C})$ (see  the discussion in section 2.1 of \cite{categor}, or Lemma 2.8 in \cite{ASnew}). The composition of the Floer product and the inverse of this isomorphism gives the desired product.

  Moreover, as in the case of symplectic cohomology, positivity of $H$ ensures the existence of continuation chain maps $\kappa^m_n: CF_m(\overline{A}, \overline{B}) \to CF_n(\overline{A}, \overline{B})$ for all $n>m$, so that $CF_m(\overline{A}, \overline{B})$ and $HF_m(\overline{A}, \overline{B})$ form a direct system (same as in the case of the Liouville domains, considered above).

 Again, we define the wrapped Floer cohomology of $\overline{A}$ and $\overline{B}$ as the direct limit of the cohomology system, $HW(\overline{A}, \overline{B})=\varinjlim HF_m(\overline{A}, \overline{B})$.

Just as in the case of admissible Lagrangians in Lioville manifolds, in addition to thimbles, we may allow either or both $L_0$ and $L_1$ to be compact Lagrangians and $HW(L_0, L_1)$ is defined in the same way.

The same properties as before hold for the wrapped Floer cohomology of thimbles, and make compact Lagrangians and Lefschetz thimbles of $E$ into objects of a Donaldson-Fukaya category.

We also note that in our circumstances of contractible thimbles the Stiefel-Whitney classes of the Lagrangians are manifestly zero, and the relative Chern class lives in $H_2 (W)$, which is zero if dim$W > 6$, 
 so in that case Floer cohomology is defined over $\mathbb{Z}$ and $\mathbb{Z}$-graded, although since our proof is based only on computing the ranks of various Floer cohomology groups, this is of marginal importance. In fact since we use the Seidel long exact sequence from \cite{S2}, we should work with $\mathbb{Z}/2$ coefficients and ungraded cohomology groups.

\begin{rem} Since a Lefschetz fibration can be made into a Lioville domain by corner smoothing, in addition to the wrapped Floer cohomology of the thimble as defined above, we can also consider the wrapped Floer cohomology of thimbles in the resulting Liouville domain. We want to comment on the relation between these two versions of the wrapped Floer cohomology.

First, we want to assume that the corner smoothings happenes in sufficiently small neighborhood of the horizontal boundary of the fibration which does not intersect either the thimbles themselves, or their images under the wrapping flow of $\hat{H}_m$ (all Lefschetz fibrations constructed using Lemma \ref{constrlem} can be arranged to have such a neigbourhood; in general one can deform a given Lefschetz fibration into one with this porperty).

Them the vanishing cycles of our thimbles are Legendrian in the smoothing. In additon, by Observation \ref{lift} one sees that  the Hamiltonian vector field $X_H$ is in the kernel of $\theta|_{\partial E}$,  so that outside of the smoothing region it is positively proportional to the Reeb vector field of the smoothing. Then there are continuation maps from $HW_l$ to $HF^{kl}$ to $HW_{k^2 l}$ for large enough $k$, which intertwine the directed systems defining the two wrapped Floer homologies, showing their equivalence. There is an additional subtlety, in that the almost complex structure used to define the wrapped Floer cohomology for thimbles is not of contact type for the smoothing; this is remedied by noting that there exists a compact subset inside the total space which contains all the holomorphic curves used in the computation of $HW$ for thimbles. Modifying the almost complex-structure outside this set to make it contact-type does not change the Floer complexes, hence the intertwining continuation maps above can be defined. The details of all these constructions (even in more generality than here) are to appear in \cite{ASnew}.   We note that nothing in the present paper relies on this equivalence of the two wrapped Floer cohomologies.

\end{rem}

\end{subsection}

\end{section}

\begin{section}{Distinguishing $W_0$ and $W_1$.}\label{sectDist}

Recall that in Section \ref{sectConstr}  we constructed, for all even $n$, two Lefschetz fibrations with total spaces $W_0$ and $W_1$ which are both diffeomorphic to $T^* S^{n+1}$ with a single subcritical handle added. In this section we prove that the wrapped Floer cohomologies of the Lefschetz thimbles are non-zero for $W_0$ and zero for $W_1$. This proves Theorem \ref{mainT}.

Briefly, the plan is to relate the computation of Floer cohomology for thimbles ``after the first wrappping" in the total space with a computation of Floer cohomology for corresponding vanishing cycles in the fiber. The resulting information about continuation map on Floer cohomlogies $HF_0(B , B) \mapsto HF_1(B ,B)$ then gives $HW(B,B)=0$.

Recall that $\tau_A$ is the Dehn twist around the sphere $A$. Our first  goal is to compute the rank of $HF(B, \tau_A(B))$ for Lagrangina spheres $A$ and $B$  in two cases - one when $A$ and $B$ are Hamiltonian isotopic, and another when they are ``sufficiently different" in the sense that will become clear shortly. These correspond to the situations with the vanishing cycles in the Lefschetz fibrations for $W_1$ and $W_2$ we have constructed in Section \ref{sectConstr}.

Of course, when $A$ is Hamiltonian isotopic to $B$ we have $HF( B, \tau_A(B))=HF( B, \tau_B(B))=HF(B, B)= H(S^{n})$ as before, and so has rank $2$.
 
 For the other case, we will use the Seidel's exact triangle in Floer cohomology  (\cite{S2}, Theorem 1). Namely, for any exact Lagrangian submanifolds $L_0$ and $L_1$ and any framed Lagrangian sphere $L$ in an Liouville domain $M$ (which for us will be a fiber of a Lefschetz fibration; this is the ``exact triangle in the fiber") we have:

 $$\to HF( L, L_1) \otimes  HF(L_0, L) \overset{m}{\to} HF(\tau_L (L_0), L_1) \to HF(L_,  L_1)\to$$
 
 It is part of Seidel's theorem that the map $m$ is the composition of isomorphisms $ HF(L_0, L)=HF(\tau_L (L_0), \tau_L (L))= HF(\tau_L (L_0), L)$ (first one is naturality of Floer cohomology, second follows from the fact that $\tau_L =L$ for all $L$), with the usual product in Floer cohomology.

 When $L=A, L_0=\tau^{-1}_A (B), L_1=B$ we get

 $$\to HF( A, B) \otimes  HF(\tau^{-1}_A (B), A)\overset{m}{\to}HF(B, B) \overset{a}{\to} HF(\tau^{-1}_A (B), B)\to$$
 
The chain of isomorphisms $ HF(L_0, L)=HF(\tau_L (L_0), \tau_L (L))= HF(\tau_L (L_0), L)$ from before becomes $ HF(\tau^{-1}_A (B), A)=HF(\tau_L \tau^{-1}_A (B), \tau_A (A))= HF(\tau_A \tau^{-1}_A (B), A)=HF(B, A)$. Using this and the isomorphism $HF(\tau^{-1}_A (B), B)=HF (B, \tau_A (B))$ we rewrite the exact triangle as

 $$\to HF( A, B) \otimes  HF(B, A)\overset{\mu}{\to}HF(B, B) \overset{\alpha}{\to} HF(B, \tau_A (B))\to$$
 
 Observe that now $\mu$ is precisely the Floer product. 
 
The reason this can be used to compute the rank of $HF(B, \tau_A (B))$ is encapsulated in the following elementary observation.

\begin{obs}\label{obsRanks} If $$\to K \overset{F}{\to} L \to M\to$$ is an exact triangle with $\operatorname{rank}K=k$, $\operatorname{rank}L=L$ and $\operatorname{rank} \operatorname{Im}F=f$, then 
$\operatorname{rank}M=k+l-2f.$

\end{obs} 
 
 Hence all we need to know are the ranks of $ HF(B, A)$, $HF(A, B)$ and the map $\mu$.

 We will have several occasions to use the following lemma relating Floer cohomology in the total space of a Lefschetz fibration and the one in thefiber.
 
  \begin{lemma} \label{lemHFcompare} Let $\pi:E \to D$ be a Lefschetz fibration. Let $\sigma$  be vanishing path and $\tau$ either a vanishing path for a different critical point or a matching path between a pair of critical points both different from the critical point of $\sigma$.  Assume $\sigma$ and $\tau$ intersect transversely exactly once at point $p$; let $S$ and $T$ be the vanishing cycles of $\sigma$ and $\tau$, respectively, in the fiber $F_p$ and let $\overline{S}$ and $\overline{T}$ be the corresponding Lefschetz thimbles. Then we have $HF(S, T)=HF(\overline{S}, \overline{T})$, where one cohomology is computed in the fiber and another in the total space.
  
  \end{lemma}
  
  \begin{proof}

The Floer chain complex computing $HF(\overline{S},\overline{T})$ has the same generators as the Floer chain complex computing  $HF(S,T)$, lying in the fiber $F_p$ over the unique intersection point $p$.  We can pick a family $J_t$ of almost complex structures on $E$ which makes the fibration map $\pi$ holomorphic, and such that their restriction to the fiber $F_p$ is regular in the sense of Floer cohomology computations in $F_p$ - that is all the pseudoholomorphic strips are regular. We claim that this family $J_t$ is then regular in the total space as well. Since any other holomorphic strips are prohibited (as they would have to project to holomorphic strips in the base), this would imply that the two Floer cohomologies - $Hom_{FS}(\overline{A},\overline{B})$ and $HF(A,B)$ are equal, and so $\operatorname{rank} Hom_{FS}(\overline{A},\overline{B})=2$.

 We want to see that $J_t$ is regular for holomorphic strips in the total space. The argument we will use has appeared  in \cite[Proposition 4.1]{MS}. We provide it here for completeness. Let $u:\mathbb{R}\times [0,1] \mapsto F_p \subset E$ be such a strip. The linearization of the operator associated to it as a map to $E$ is a Fredholm operator $D_u:\mathcal{H}_1 \to \mathcal{H}_0$. If we consider the same $u$ as a map to $F_p$ the linearization  is the restriction of $D_u$ to subspaces  $\overline{\mathcal{H}}_1 \to \overline{\mathcal{H}}_0$. By projecting to the base, we identify the quotient space $\mathcal{H}_1/\overline{\mathcal{H}}_1$ with $W^{1,q}$ (for some $q\geq2$)  space of functions  $\xi: \mathbb{R} \times [0,1] \rightarrow \mathbb{C}$,  satisfying boundary conditions $\xi(\mathbb{R} \times \{0\}) \subset \sigma_0\mathbb{R}$, $\xi(\mathbb{R} \times \{1\}) \in {\sigma_1}\mathbb{R}$. Here, $\sigma_0\mathbb{R} \subset \mathbb{C}$ and $\sigma_1 \mathbb{R} \subset \mathbb{C}$ are the tangent spaces of $\sigma$ and $\tau$ at $p$, hence transverse by assumption. Similarly, $\mathcal{H}_0/\overline{\mathcal{H}}_0$ can be identified with the space of all $L^q$ functions $\mathbb{R} \times [0,1] \rightarrow \mathbb{C}$, and the quotient map induced by $D_u$ is the standard Cauchy-Riemann operator $\bar\partial$, which is invertible. This implies that regularity in $F_p$ and in $E$ are equivalent, as claimed.

  \end{proof}

 We summarize the first computation in the following lemma.
 
 \begin{lemma} \label{lemrank0} Let $W_1$ be the total space of the Lefschetz fibration constructed in Section \ref{sectConstr}. Then for the two vanishing cycles $A$ and $B$ in the fiber, we have  $\operatorname{rank} HF(B, \tau_A (B))=\operatorname{rank} HF( \tau_A (B), B)= 4$. 
  
 \end{lemma}

 \begin{proof}
  Note that we still have $HF(B,B)= H(S^{n})$ has rank $2$.

The spheres $A$ and $B$ intersect transversely at the two critical points of $\rho$ which have the same grading mod 2. This can be seen, for example,  since gradings are unchanged by totally-real isotopy and $B$ is (non-Lagrangian!) totally-real isotopic to $A$, so the gradings on the intersection points are the same as gradings of intersections of $A$ and a perturbation of $A$, which by Morse-Bott techniques going back to the same already-cited computation of Floer in \cite{Floer2}, differ by the dimension of the sphere $A$, an even number.
Hence the Floer differential vanishes, and rank $HF(A,B)$ is equal to 2 even for the case of $W_1$.
 
 On  the other hand, the vanishing cycles $A$ and $B$ are not Hamiltonian isotopic. In fact, they are not isomorphic as objects of the Donaldson-Fukaya category of the fiber. To see this consider the Lefschetz thimble $L$ for the critical point over $0$ of the auxiliary fibration $\rho$ and the vanishing path going straight down. We see that in the total space of $\rho$, i.e the fiber of $\pi$,  $A$ and $L$ are disjoint and so $HF(A, L)=0$. On the other hand, the vanishing path for $L$ intersects the matching path for $B$ exactly once, so by Lemma \ref{lemHFcompare} applied to $L$ and $B$ we see that $HF(B, L)$ is the same as the Floer cohomology of the corresponding vanishing cycles (of $\rho$), that is of the vanishing sphere with itself.  This is again $H(S^{n-1})$, and so is not zero. 
  
 Correspondingly, for these non-isomorphic vanishing cycles the pair of pants product still hits the fundamental class in $HF(B,B)$ (by Poincare duality in Floer theory, \cite{Book}[Sections 8c and 12e]). We want to see that the product  misses the identity, which in turn by  Observation \ref{obsRanks} forces the group $HF(B, \tau_A (B))$ to be of rank $4+2-2=4$.  Of course,  $\operatorname{rank} HF(B, \tau_A (B))=\operatorname{rank} HF( \tau_A (B), B)$ by Poincare duality in Floer cohmology.
 
 To see that $\mu$ misses the identity, suppose there is $c \in HF( A, B) \otimes  HF(B, A)$ such that $\mu(c)= Id \in HF(B,B)$. Then for a non-zero element $d \in HF(B, L)$, the composition of Floer products going from $HF(B, L) \otimes HF( A, B) \otimes  HF(B, A)$ to $HF(B, L)$ on the one hand takes $d \otimes c$ to $\mu(d, \operatorname{Id})=d$, and on the other hand factors through the group $HF(A, L)=0$.  As $d \neq 0$, this is a contradiction.

 \end{proof}
 
Next we would like to understand the Floer cohomology of thimbles in the total spaces of $W_0$ and $W_1$. In both cases we  denote the thimbles of the Lefschetz fibration by $\overline{A}$ and $\overline{B}$ and the vanishing paths they lie over by $\alpha$ and $\beta$ respectively, with critical values $\alpha(1)=a$ and $\beta(1)=b$; further, in notation of Section \ref{sectWHF},  we write $\overline{B_1} =\phi^1 (\overline{B})$; that is $\overline{B_1}$ is the result of wrapping $\overline{B}$ once around the base, the thimble lying over a wrapped path $\beta_1$.

 \begin{lemma} \label{lemRank} In both $W_0$ and $W_1$  we have  $\operatorname{rank} HF(\overline{B_1}, \overline{B})= 3$.
 \end{lemma}
 
 \begin{proof}

We view $\overline{A}$ and $\overline{B}$ and as objects of the derived Fukaya-Seidel category of the Lefschetz fibration $W_i \to D^2$ (\cite{Book}[Section 18f]; Seidel calls it the Fukaya category of a Lefschetz fibraion, denoted $\mathcal{F}(\pi)$).   
 
The results of \cite{VCM} (see also \cite{Book}[Proposition 18.23]) imply that $\overline{B_1}$ is isomorphic to the cone of the evaluation map $ev: Hom(\overline{A},\overline{B})\otimes \overline{A} \mapsto \overline{B}$.  

Taking the corresponding exact triangle and taking the long exact sequence corresponding to applying the functor  $Hom(\cdot,\overline{B})$ to it we get in cohomology

 $$ \to Hom_{FS}(\overline{B_1}, \overline{B}) \to Hom_{FS}(\overline{B},\overline{B}) \overset{ev^*}{ \to} Hom_{FS}(\overline{A},\overline{B}) ^*\otimes Hom_{FS}(\overline{A},\overline{B}) \to $$ 
 
Here for thimbles $X$ and $Y$ we write $Hom_{FS}(X,Y)=HF(X,Y)$ for morphisms in cohomology level Fukaya-Seidel category, that is the Floer cohomologies. We keep this notation to better distinguish Floer cohomology in the Fiber (still denoted $HF$) and the total space (denoted $Hom_{FS}$).
 
We now want to use Observation \ref{obsRanks} to compute $\operatorname{rank} Hom_{FS}(\overline{B_1}, \overline{B})$.

   $Hom_{FS}(\overline{B},\overline{B})$ has rank one and is generated by the identity. This is a general fact about Lefschetz thimbles - a vanishing path can be isotoped to intersect itself only at the critical point, which changes the thimble by a Hamiltonian isotopy (\cite{Aur}, Lemma 3.2), so the resulting cohomology group has rank 1.

 Now by Lemma \ref{lemHFcompare} applied to the thimbles $\overline{A}$ and $\overline{B}$ we get $Hom_{FS}(\overline{A},\overline{B})=HF(A,B)$, and so $\operatorname{rank} Hom_{FS}(\overline{A},\overline{B})=2$.

 The map $ev^{*}$ is non-zero and maps $Id \mapsto \Sigma_{\alpha \in  Hom_{FS}(\overline{A},\overline{B})} \alpha \bigotimes \alpha^{*}$. By the Observation \ref{obsRanks} this means the $\operatorname{rank }Hom_{FS} (\overline{B_1}, \overline{B})=1+4-2=3$, as wanted.

 \end{proof}

 We now use this information to prove our main theorem.
 
 \begin{proof}[Proof of Theorem \ref{mainT}]

We would like to understand the Floer cohomology $Hom_{FS} (\overline{B_1}, \overline{B})$. The thimble $B_1$ lies over $\beta_1$ and $B$ lies over $\beta$, so the chain complex computing $Hom_{FS} (\overline{B_1}, \overline{B})$ contains the generators lying over $\beta \cap \beta_1$, that  is one generator $e$ corresponding to the critical point of the main fibration, where the thimbles $\overline{B}$ and $\overline{B_1}$ meet, and the generators lying over the other intersection point $x$, which correspond to the generators of the complex computing the group $HF(\tau_A (B), B)$.

We would like to use the projection $\pi$ again. However, the family of almost complex structures $J_t$ that makes $\pi$ holomorphic may not be regular for the holomorphic strips not contained in the fiber. Nonetheless, we could still use it to get some information. In particular, for such a $J_t$ there are no holomorphic strips from $e$ to generators in $F_x$ (their projections would be holomorphic strips going from $e$ to $x$, which don't exist). Then by Gromov compactness, we conclude that for a sufficiently small regular perturbation of the family $J_t$ there are no such holomorphic curves either (else they would subconverge to a curve for $J_t$ itself). This means that those generators of the Floer complex of $\overline{B}$ and $\overline{B}_1$  lying over $F_x$ form a subcomplex (it is important here that we are doing Floer \emph{cohomology}, so terms in the differential of a generator $g$ are given by holomorphic strips going \emph{to} $g$.)

Further, as in the proof of Lemma \ref{lemHFcompare}, if the restriction of the original family $J_t$ is regular for the computation of $HF(\tau_A (B), B)$ in the fiber $F_x$, it is also regular for those holomorphic strips viewed as lying in the total space $W_i$. This means that for small perturbation $\hat{J_t}$  there is one to one correspondence between the  $\hat{J_t}$-holomorphic strips contained in $F_x$  and the $\hat{J_t}$-holomorphic strips contained in $F_x$. This in turn means that the subcomplex of generators lying in $F_x$ is isomorphic to the complex computing the Floer cohomology $HF(\tau_A (B), B)$ in the fiber.

Basic homological algebra then says that the cohomology $Hom_{FS} (\overline{B_1}, \overline{B})$ is computed by the complex $\mathbb{Z}_2 \overset{\epsilon}{\to} HF(\tau_A (B), B)$. We note that in this situation there is always a map $i:\mathbb{Z}_2 \to Hom_{FS} (\overline{B_1}, \overline{B})$ induced by the inclusion of $\mathbb{Z}_2$ to the first term of the complex. 

We have the following lemma.

\begin{lemma} The  map $i:\mathbb{Z}_2 \to Hom_{FS} (\overline{B_1}, \overline{B})$ above coincides with the continuation map $\mathbb{Z}_2= Hom_{FS} (\overline{B}, \overline{B}) \to Hom_{FS} (\overline{B_1}, \overline{B})$
\end{lemma}

\begin{proof} This is Proposition 4.2 in \cite{MS}. The outline is as follows. Writing down the continuation map equation, and assuming that the family $J_{s,t}$ of almost complex structures used in it makes $\pi$-holomorphic, one gets that  the projection of any solution is a solution of a ``continuation map" type equation in $\mathbb{C}$. An energy estimate (in $\mathbb{C}$) shows that the only solutions of that are constant (at $b$), so that solutions to the original continuation map equation are contained in in the fiber $F_b$. By a Gromov compactness argument as before, one concludes that the same is true for a generic small perturbation of $J_{s,t}$, so that solutions of continuation map equation are contained  $F_b$ even for the perturbed family. In $F_b$ the continuation map equation becomes the Cauchy-Riemann equation and has unique regular solution - namely the constant one. Hence the continuation map is given on the chain level by including the generator $e$ into the chain complex, as wanted. We refer the reader to \cite[Proposition 4.2]{MS} for details.  
\end{proof}

We now concentrate on the more ``exotic" manifold $W_1$.
 
Since in $W_1$ the rank of $HF(\tau_A (B), B)$ is 4 by Lemma \ref{lemrank0} and  $\operatorname{rank} HF(\overline{B}, \overline{B_1})= 3$ by Lemma \ref{lemRank}, we conclude that the map $\epsilon: \mathbb{Z}_2 \to HF(\tau_A (B), B)$ is non-zero, and the Floer cohomology $Hom_{FS} (\overline{B_1}, \overline{B})$ is $HF(\tau_A (B), B)/ \operatorname{Im} \epsilon$, so that the map $i$ is the zero map.

However, looking back to the definition of the product structure on $HW(\overline{B}, \overline{B})$,  we see that the unit is given by the image of the generator  $u \in Hom_{FS}(\overline{B},\overline{B})$ under continuation maps. Since $i(u)=0$,  the unit in $HW(\overline{B}, \overline{B})$ is zero, which is only possible if $HW(\overline{B},\overline{B})=0$.  
Since $HW(\overline{A}, \overline{B})$ is a module over $HW(\overline{B},\overline{B})$ it also vanishes. A symmetric argument implies that  $HW(\overline{A},\overline{A})=0$.

\begin{rem}We note that this behavior of wrapped Floer cohomology is in sharp contrast to the one in standard cotangent bundles. There  for  $F_i$ the cotangent fiber at the point $p_i$ the wrapped Floer cohomology $HW(F_1,  F_2)$ is the homology of the path space from $p_1$ to $p_2$ (see Theorem 3.2 in \cite{categor}).  In fact by analogy with the result of Cieliebak  that states that subcritical handle attachment does not change the symplectic cohomology (\cite{cieliebak}), we expect that in the case of subcritical  handle attachment  the functor constructed by Abouzaid and Seidel in \cite{AS}  is a full embedding, in which case the wrapped Floer cohomologies in $W_0$ should coincide with those in $T^* S^{n+1}$ from which it is obtained (see \cite{AS2}[Property 2.5] for a related result). 
\end{rem}

Meanwhile, we conclude that the manifold $W_1$ does not contain any closed exact Lagrangian submanifold. To see this, observe that, for such a Lagrangian $L$ we would have on one hand, by Floer's original result cited before, $HF(L, L)=H^*(L)$,  which is nonzero. On the other hand,  $HF(L,L)=HW(L,L)$ since wrapping does not affect closed Lagrangian submanifolds, but by Theorem 4 of \cite{FSS} there is a spectral sequence converging to $HF(L,L)$ with the first page $E_1^{jk}=(HF(\Delta^!_j, L)\bigotimes HF(L,\Delta_j ))^{j+k}$ for a basis of thimbles $\Delta$ and dual thimbles $\Delta^!$. However, $HF(L, \Delta_j)=HW(L, \Delta_j)$, again because $L$ is closed, and the later group vanishes since it is a module over $HW(\Delta_j, \Delta_j)=0$. Hence, the above results imply that this first page vanishes, a contradiction.

As $W_0$ contains the exact Lagrangian sphere inherited from the zero-section of $T^* S^{n+1}$, we conclude that $W_0$ and $W_1$ are not Liouville isomorphic.

We also note that were the wrapped Floer cohomology $HW(\overline{A}, \overline{A})$ to vanish in the case of $W_0$, then by symmetry so would $HW(\overline{B}, \overline{B})$, and we could repeat the above argument. So the fact that there is an exact Lagrangian sphere in $W_0$ implies that these groups are non-zero.

This completes the proof.

\end{proof}

\begin{rem} After the work on this paper has been completed, Abouzaid and Seidel in \cite[Property 2.3]{AS2} have proved that vanishing of Wrapped Floer cohomology for all thimbles implies vanishing of symplectic cohomology, which in turn by a result of Viterbo implies that the Liouville domain contains no compact exact Lagrangians (see \cite{Biased}[Theorem 5.1]). This gives an alternative way to prove that $W_1$ has no compact exact Lagrangian submanifolds. 
\end{rem}

\end{section}

\begin{section}{Extensions and further research.}

 After the work on this paper has been completed, the methods used here were extended first in joint work with Seidel 
 \cite{MS}, and then quite radically by  Abouzaid and Seidel in \cite{AS2}, where they have proved the following vast generalization.
 
 \begin{tm}{\cite[Theorem 1.1]{AS2}} For any Liouville domain $W$ of dimension $\geq 6$ there exists a Liouville domain $W^{\prime}$ such that

\begin{itemize}
\item $W$ and $W^{\prime}$ are diffeomorphic
\item The symplectic cohomology and wrapped Fukaya category of $W^{\prime}$ are zero.
\end{itemize}

 \end{tm}

 Moreover, they obtain for any   Liouville domain $W$ a sequence of domains $W_k$, all of which are almost diffeomorphic to $W$, and such that $W_k$ is not Liouville-isomorphic  to $W_l$ for $k \neq l$.

 These results show the vastness of the world of ``exotic" symplectic manifolds in general, and of ``empty" ones (those with vanishing Floer--theoretic invariants) in particular.

Not much is known about distinguishing ``empty" symplectic manifolds. One positive result is the work of Richard Harris, who, building on methods of this paper,  has shown in \cite[Theorem 1.3]{harris}, that there exist for any $n \geq 1$ a manifold $W$ (diffeomorphic to $T^*S^3$ with $n$ 2-handles attached), and exact symplectic forms $\omega_1, \ldots \omega_{n+1}$  on it such that, with respect to each $\omega_i$ the manifold $(W, \omega_i$) is Liouville and contains no exact Lagrangian submanifolds, but such that there exists no diffeomorphism $\phi$ of $W$ such that $\phi^* \omega_i = \omega_j$ for $i\neq j$.  These $(W, \omega_i)$ all have vanishing wrapped Fukaya categories and vanishing symplectic cohomology.
 
 Harris's method relies on studying a non-exact deformation of the symplectic structures on the fibers of relevant Lefschetz fibrations, and as such is not directly applicable to higher-dimensional examples.  I hope to develop more algebraic versions of Harris's method in the future.

\end{section}

 \begin{section}*{Acknowledgments}

 This paper is based on the author's Ph.D. thesis written under supervision of Denis Auroux. I would also like to thank Paul Seidel for inspiring this research problem. This area of research owes much to his ideas, as should be manifest in the present paper as well.

\end{section}

\end{document}